\documentclass[11pt,reqno]{amsart}

\pdfoutput=1
\numberwithin{equation}{section} \oddsidemargin=-.0cm
\usepackage{amssymb,amsmath,amsfonts}

\usepackage{newtxtext}
\usepackage[scaled=0.90]{helvet} 
\usepackage{array, booktabs}  
\usepackage{color}
\usepackage{sidecap}

\usepackage[T1]{fontenc}
\usepackage{amscd,mathtools}
\usepackage{stmaryrd}
\usepackage{dsfont}
\usepackage{needspace} 
\usepackage{cases} 
\sidecaptionvpos{figure}{c}
\usepackage{graphicx, graphics}
\usepackage{mathrsfs}

\usepackage[sort]{cite}
  \usepackage{paralist}
  \usepackage{graphics} 
  \usepackage{epsfig} 
\usepackage{graphicx}  
\usepackage{epstopdf}
\usepackage{bbm} 

\usepackage[colorlinks=true, pdfstartview=FitV, linkcolor=blue, 
                     citecolor=blue, urlcolor=blue]{hyperref}  

\newcommand{\ep}{\varepsilon}

\graphicspath{{./Figures/}}

\def\Lip{{\rm Lip}}

\def\cH{\mathcal H}
\def\Forall{\text{ } \forall \:}
\def\d{\mathrm{d}}

\newtheorem{thm}{Theorem}[section]
\newtheorem{lem}{Lemma}[section]
\newtheorem{defi}{Definition}[section]

\newtheorem{prop}[thm]{Proposition}

\newtheorem{rem}{Remark}[section]
\newtheorem{cor}{Corollary}[section]
\def\bt{\begin{thm}}
\def\et{\end{thm}}
\def\bl{\begin{lem}}
\def\el{\end{lem}}
\def\bd{\begin{defi}}
\def\ed{\end{defi}}
\def\bc{\begin{cor}}
\def\ec{\end{cor}}
\def\bp{\begin{proof}}
\def\ep{\end{proof}}
\def\br{\begin{rem}}
\def\er{\end{rem}}

\def\bprop{\begin{prop}}
\def\eprop{\end{prop}}

\def\Forall{\text{ } \forall \:}
\def\d{\, \mathrm{d}}

\def\be{\begin{equation}}
\def\ee{\end{equation}}
\def\bes{\begin{equation*}}
\def\ees{\end{equation*}}
\def\bea{\begin{equation} \begin{aligned}}
\def\eea{\end{aligned} \end{equation}}
\def\beas{\begin{equation*} \begin{aligned}}
\def\eeas{\end{aligned} \end{equation*}}

\def\bpp{\begin{proof}}
\def\epp{\end{proof}}

\def\bi{\begin{itemize}}
\def\ei{\end{itemize}}
\def\ben{\begin{enumerate}}
\def\een{\end{enumerate}}

\date{July 22, 2020}
\title[Optimal management of harvested population at the edge of extinction]{Optimal management of harvested population at the edge of extinction}

\author[Micka\"el D. Chekroun]{Micka\"el D. Chekroun}
\address[MDC]{Department of Earth and Planetary Sciences, Weizmann Institute, Rehovot 76100, Israel;
Department of Atmospheric and Oceanic Sciences, University of California, Los Angeles, CA 90095-1565, USA}
\email{mchekroun@atmos.ucla.edu}

\author[Honghu Liu]{Honghu Liu}
\address[HL]{Department of Mathematics, Virginia Tech, Blacksburg, Virginia 24061, USA}
\email{hhliu@vt.edu}

\begin{document}

\begin{abstract} 
Optimal control of  harvested population at the edge of extinction in an unprotected area, is considered. The underlying population dynamics is governed by
a Kolmogorov-Petrovsky-Piskunov equation with a harvesting term and space-dependent coefficients while the control consists of transporting individuals from a natural reserve. The nonlinear optimal control problem is approximated by means of a Galerkin scheme.  Convergence result about the optimal controlled solutions and error estimates between the corresponding optimal controls,  are derived.  For certain parameter regimes, nearly optimal solutions are calculated from a simple logistic ordinary differential equation (ODE) with a harvesting term, obtained as a Galerkin approximation of the original partial differential equation (PDE) model. A critical  allowable fraction $\underline{\alpha}$ of the reserve's population  is inferred from  the reduced  logistic ODE with a harvesting term. This estimate obtained from the reduced model allows us to distinguish sharply between survival and extinction for the full PDE itself, and thus to declare whether a control strategy leads to success or failure for the corresponding rescue operation while ensuring survival in the reserve's population. In dynamical terms, this result illustrates that although continuous dependence on the forcing may hold on finite-time intervals, a high sensitivity in the system's response may occur in the asymptotic time.  We believe that this work, by its generality, establishes  bridges interesting to explore between optimal control problems of ODEs with a harvesting term  and their PDE counterpart.
\end{abstract}

\maketitle
\tableofcontents

\section{Introduction}
Overexploitation has led to the extinction of many
species~\cite{ba}. Traditionally, models of ordinary differential
equations (ODEs) or difference equations have been used to estimate
the maximum sustainable yields from populations and to perform
quantitative analysis of harvesting policies and management
strategies~\cite{getz}. Ignoring age or stage structures as well as
delay mechanisms, which will not be treated by the present paper,
the ODE models are generally of the type 
\be \label{eqi1} 
\frac{\d
U}{\d t}=F(U)-Y(U), 
\ee 
where $U$ is the population
biomass at time $t$, $F(U)$ is the growth function, and $Y(U)$
corresponds to the harvest function. In these models, the most
commonly used growth function is logistic, with $F(U)=U(\mu-\nu U)$
\cite{bedmay,schaefer1954some}, where $\mu>0$ is the
intrinsic growth rate of the population and $\nu>0$ models its
susceptibility to crowding effects.

Different harvesting strategies $Y(U)$ have been considered in the
literature and are used in practical resource management. A very
common one is the \emph{constant-yield harvesting} strategy, where a
constant number of individuals are removed per unit of time:
$Y(U)=\delta$, with $\delta$ a positive constant. This harvesting
function naturally appears when a quota is set on the
harvesters~\cite{robinson1999towards,stephens2002sustainable}. Another
frequently used harvesting strategy is the \emph{proportional
harvesting} strategy (also called \emph{constant-effort
harvesting}), where a constant proportion of the population is
removed. It leads to a harvesting function of the type $Y(U)=\delta
U$.

Reaction-diffusion equations have also been used extensively in modeling the spatiotemporal behavior of a species of organism \cite{shigesada1997biological,cantrell2004spatial,murray2001mathematical,murray2007mathematical,okubo2013diffusion},
and harvesting effects have been analyzed within this class of models; see e.g.~\cite{Oruganti2002,RC07,RC10}.
Within this global picture, we focus on a particular class of models, namely the class of diffusive logistic equation in a
heterogenous environment \cite{skt,cantrell1989spatial,shigesada1997biological,bh,bhr1,bhr2}. 

More specifically, the harvesting model we consider, described in Sec.~\ref{Sec_KPPH} below, is the {\bf  Kolmogorov-Petrovsky-Piskunov harvesting (KPPH)} model (see \eqref{Eq_KPP} below) in which the local growth rate of the population and the limit effects to crowding are spatial-dependent coefficients. The harvest function is a quasi-constant-yield harvesting term depending on a threshold below which harvesting is progressively abandoned. This model is analyzed in \cite{chekroun2006models,RC07} and the main results about its steady state analysis and its dynamical behavior, are recalled in Sec.~\ref{Sec_background2} below for the reader's convenience. In particular, the asymptotic behavior of the KPPH model was characterized as a function of
the harvesting  intensity $\delta$; see \eqref{Eq_KPP}.  It was proved in \cite{RC07}  that if $\delta$ is smaller than a critical value $\delta^*$, the
population density converges to a ``significant'' state, which is everywhere above a small threshold. On the other hand, it was shown in  \cite{RC07} that if $\delta$ is
larger than another critical value $\delta_2$, which is slightly above $\delta^*$ in practice, the population density eventually settles down to a
``remnant'' state, everywhere below the same small threshold. The population can be considered as extinct in this case.  

Within this context, we assume that given an unprotected area $\Omega$ in which the population (of let's say some mammals) evolves according to a KPPH equation, is under threat of extinction due to some illegal hunting that superimposes to some allowed  harvesting, causing thereof a rise in the harvesting intensity such that  $\delta$ ends to jump above $\delta_2$. Given a natural reserve $\widehat{\Omega}$  in which the same species evolves, we address in this article the problem of saving the population under threat  by releasing in a controlled fashion new individuals from the reserve, into the unprotected area, while avoiding to exert too much pressure on the reserve that would be detrimental on it.  Our control strategy is seriously constrained in time as the operation is assumed to start when a (large) fraction of the original population has been decimated and its course to extinction is, thus, well engaged.  

The resulting optimal control problem differs fundamentally from other optimal control problems concerned with harvesting models that appeared in previous studies which dealt mainly with the search for optimal harvesting strategies; see e.g.~\cite{neub,kurata2008optimal,finotti2012optimal,Brock2014optimal}. Instead, we  aim at designing optimal rescue strategies while taking into account plausible factors that may arise in real-life problems, from a reserve management perspective.  

The originality of this work does not only rely  on its problem formulation, but also on its proposed solution. 
We indeed address the obtention of optimal solutions by means of rigorous Galerkin approximations of the underlying KPPH models. To do so, we rely on 
the recent mathematical framework introduced in \cite{CKL17} which allows for deriving convergence results and error estimates from Galerkin approximations of a broad class of nonlinear optimal control problems in infinite dimension.  Of course, the idea of designing approximations of optimal solutions from finite-dimensional approximations is not new, and a great deal of works have addressed this question in various contexts \cite{Ded10,GK11,Hinze_al05,IK08,TV09} as well as based on various basis functions, possibly empirical \cite{Franke12,HK98,HK00,ravindran2000reduced}.  However, rigorous convergence analysis from finite-dimensional Galerkin approximations do not seem to have been much addressed for the optimal control of nonlinear problems, and in that sense \cite{CKL17} provides some useful elements. 

In that respect, we show that the framework of \cite{CKL17} allows us to derive error estimates between the optimal control of the KPPH model considered in this article and the optimal control built from any of its Galerkin approximations.  These error estimates are supplemented by convergence of the corresponding 
controlled solutions towards the optimally controlled solution of the full partial differential equation (PDE) model. These convergence and error estimates results are summarized in  Theorem \ref{Thm_cve_KPP} below which is proved in Appendix  \ref{Sect_cond_verification} based on the theoretical apparatus from \cite{CKL17}, recalled in Appendix~\ref{Sect_theory}. 

The finite-dimensional approximations are here simply obtained from projection onto the eigenmodes of a natural underlying spectral problem (see \eqref{Eq_spectral_pb} below). A standard Pontryagin Maximum Principle (PMP) approach \cite{PBGM64,Kirk12} is then used to find extremals of the finite-dimensional optimal control problems \cite{bonnard2003singular},  which turn out to approximate,  in our case, the optimal solutions of the full problem; see  Sec.~\ref{Section_PMP}. 

To illustrate our theoretical framework and to favor reproducibility of the results, the numerical results of Sec.~\ref{Sec_numerical_results} are presented within the context of one-dimensional environments and for various level of fragmentation of the population's habitat.  For the parameter regime considered, we show 
that from a one-dimensional, ordinary differential equation (ODE) approximation, nearly optimal solutions  can be derived. The corresponding Galerkin approximation of the KPPH model, reducing to a simple logistic ODE with a harvesting term\footnote{Whose solution contains most of the PDE solution's energy.}, allows us in turn to provide great insights about the optimal control problem of the PDE model. The main result, from an ecological perspective, is indeed expressed in terms of a fraction $\alpha$ of the reserve that is allowed for exploitation.  
This fraction, related to e.g.~management policy of the reserve, is caping the amount of individuals transported from the reserve into the unprotected area. 
We show that a critical  fraction $\underline{\alpha}$ can  be inferred from the reduced  logistic ODE with a harvesting term. This estimate obtained from the reduced equation is indeed  shown to be very useful as it allows us to distinguish sharply between survival and extinction for the full PDE itself, and thus to declare whether a control strategy leads to success or failure for the corresponding rescue operation  while ensuring survival in the reserve's population.  In terms of optimal control, this result translates into the existence of  controls very close to each other, one of which leading to a significant survival of the population while the other leading to its extinction. In dynamical terms, this result illustrates that although continuous dependence on the forcing may hold on finite-time intervals, a high sensitivity in the system's response may occur in the asymptotic time.

\section{Optimal control of harvested population}

\subsection{The KPP model with harvesting term in heterogeneous environment}\label{Sec_KPPH}
The partial differential equation (PDE) and boundary conditions underlying the optimal control problem considered hereafter are described as follows. 
Let $\Omega$ be a smooth bounded and connected domain of $\mathbb{R}^d$ ($d\geq 1$). We consider
\begin{subnumcases}{\label{Eq_KPP}}
\; \partial_t y = D \nabla^2 y + \mu(x) y - \nu(x) y^2 - \delta \rho_{\epsilon}(y), \quad (t,x)\in (0,\infty) \times \Omega, \label{Eq_KPP1}\\
 \;\frac{\partial y}{\partial {\boldsymbol{n}}}= 0, \quad (t,x)\in [0,\infty) \times \partial \Omega. \label{Eq_KPP2}
\end{subnumcases}
Here $\boldsymbol{n}$ denotes the outward unit normal to  the boundary $\partial \Omega$.

This equation differs from the classical Fisher equation \cite{fisher} (also known as the Kolmogorov-Petrovsky-Piskunov (KPP) equation \cite{KPP}), by its spatially-dependent coefficients, $\mu(x)$ and $\nu(x)$, as well as its harvesting term $\delta \rho_{\epsilon}(y)$. When $\delta=0$ in \eqref{Eq_KPP}, the model reduces to the Shigesada-Kawasaki-Teramoto 
model described in \cite{skt}; see also \cite{bh,bhr1,bhr2,shigesada1997biological}. Such a problem fits with general species assessment and management problematics considered  for instance  in \cite{ba,bedmay,getz,neub,RC10, schaefer,fahrig2003effects,kurata2008optimal}. Hereafter, Eq.~\eqref{Eq_KPP1} will be referred to as the Kolmogorov-Petrovsky-Piskunov harvesting (KPPH) equation.  

The unknown function $y = y(t, x)$ denotes the population density at time $t$ and space position $x$.
The coefficient $\mu$ represents the intrinsic growth rate of the population,  which is assume to be a measurable function of $x$ in $L^\infty(\Omega)$.
 The spatial dependence of $\mu$ is introduced to account for the possible impact of environmental heterogeneity \cite{RC07}.
The function $\nu(x)>0$ (also in $L^\infty(\Omega)$) represents the susceptibility to crowding effects and is interpreted as an intraspecific competition term. 
Regions with higher values of $\mu(x)$ and lower values of $\nu(x)$
are qualified as being {\em more favorable}, while, on the other
hand, regions with lower $\mu(x)$ and higher $\nu(x)$ values are
considered as being {\em less favorable\/} or, equivalently, {\em
more hostile}.

The harvesting function $\rho_{\varepsilon}$ satisfies
\be \label{Eq_rho_conditions}
\rho_{\epsilon} \in C^1(\mathbb{R}), \quad \rho_{\epsilon}' \ge 0, \quad \rho_{\epsilon}(s) = 0,\; \; \forall s \le 0, \quad \text{and} \quad  \rho_{\epsilon}(s) = 1,\;\; \forall s \ge \epsilon,
\ee
where $\epsilon$ is a nonnegative parameter, taken to be sufficiently small in a sense made precise in Sec.~\ref{Sec_background2}.

The term $\delta \rho_{\epsilon}(y)$ with $\delta \geq 0$,
corresponds to a {\it quasi-constant-yield} harvesting term.
Indeed,  for such a harvesting
function, the yield  is constant in time whenever $y\geq
\epsilon$, while it depends on the population density when
$y<\epsilon$. Note that the
function $\rho_\epsilon$ ensures the nonnegativity of the
solutions to \eqref{Eq_KPP}; see \cite{RC07}. From a biological viewpoint,
$\varepsilon$ corresponds to a threshold below which harvesting
is progressively abandoned. Considering constant-yield harvesting
functions without this threshold value would be unrealistic since it
would eventually lead to harvest on zero-populations.

\subsection{Main results from \cite{RC07}}\label{Sec_background2}
The problem \eqref{Eq_KPP} has been analyzed in \cite{RC07}.
Using sub- and supersolution methods and the
characterization of the first eigenvalue of the linearized elliptic
operator at the trivial solution, the authors obtained existence and nonexistence results as well as
results on the number of stationary solutions; see also \cite{chekroun2006models}. 

The asymptotic behavior of the evolution equation was in particular  characterized as a function of
the harvesting  intensity $\delta$.  In \cite{RC07} it was proved that if $\delta$ is smaller than a critical value $\delta^*$, the
population density converges to a ``significant'' state, which is
everywhere above a small threshold (depending on $\epsilon$) while if $\delta$ is
larger than $\delta_2$ (another threshold that is bigger than but close to $\delta^*$), the population density $y(t,x)$ converges to a
``remnant'' state, everywhere below this small threshold. Theorem \ref{Thm_RC07_results}
below summarizes the main asymptotic and existence results from \cite{RC07}.

To formulate these results we recall some tools  used in \cite{RC07} and in the formulation of our optimal control problem below.
For this purpose, we first consider the eigenvalue problem associated with the linearization of \eqref{Eq_KPP} at the trivial solution:
\begin{subnumcases}{\label{Eq_spectral_pb}}
\;- D \nabla^2 \phi - \mu(x) \phi  = \lambda \phi, \quad  x \in \Omega, \label{Eq_spectral_pb1}\\
\;\frac{\partial \phi}{\partial {\boldsymbol{n}}}= 0, \quad  x \in \partial \Omega. \label{Eq_spectral_pb2}
\end{subnumcases}
Let  $\lambda_1$ denote  the first eigenvalue and $\phi_1$  its corresponding positive eigenfunction, unique when normalized. In particular
$\phi_1$ satisfies 
\be \label{Eq_phi1_normalization}
\phi_1 (x) > 0 \quad \forall \; x \in \Omega, \quad \text{and} \quad \|\phi_1\|_{\infty} = 1. 
\ee
The above normalization is possible since the first eigenfunction has a fixed sign in $\Omega$ as consequence of the Krein-Rutman theorem; see e.g.~\cite{amann1976fixed}.  We introduce next 
\be
\underline{\phi_1} = \min_{x \in \Omega} \phi_1,
\ee 
Note that $\underline{\phi_1} > 0$ since $\phi_1$ does not vanish on the boundary $\partial \Omega$; see \cite{amann1976fixed}. Note also that $\underline{\phi_1} \le 1$ thanks to \eqref{Eq_phi1_normalization}.

Recall from \cite[Definition 2.5]{RC07} that, a stationary solution of \eqref{Eq_KPP}, $p_{\delta}$, is called a significant solution if 
\be \label{Def_significant_soln}
\min_{x \in \Omega} p_{\delta} \ge \frac{\epsilon}{\underline{\phi_1}}.
\ee
On the other hand, a stationary solution $p_{\delta}$ is called remnant if 
\be \label{Def_remnant}
\max_{x \in \Omega} p_{\delta} < \frac{\epsilon}{\underline{\phi_1}}.
\ee

With these tools in hand,  it was derived in \cite{RC07} the following analytic formulas allowing for estimating  
the critical harvesting intensity $\delta^\ast$ leading to decline of the population towards a remnant steady state (i.e.~close to extinction):
\be \label{Eq_deltas_def}
 \quad \delta_1 = \frac{\lambda_1^2 \underline{\phi_1}}{ \overline{\nu} (1+\underline{\phi_1})^2}, \quad \delta_2 = \frac{\lambda_1^2}{4  \underline{\nu}},
\ee
see \cite[(2.15)]{RC07} with $\alpha=1$  since  $h(x)$ therein is identically equal to 1 here. 
Finally, $\overline{\nu} = \max_{x \in \Omega} \nu(x)$ and $\underline{\nu} = \min_{x \in \Omega} \nu(x)$. 

We are now in position to summarize the main theoretical results from \cite{RC07} into the following theorem.  

\bt \label{Thm_RC07_results}

{\bf Steady states \cite{RC07}.}
If $\lambda_1 < 0$, then there exists a threshold $\delta^* \ge 0$ such that 
\bi
\item[$(i)$] if $\delta \le \delta^*$,  there exists at least one positive stationary significant solution of \eqref{Eq_KPP}, whereas  
\item[$(ii)$] if $\delta > \delta^*$, there is no positive stationary significant solution of \eqref{Eq_KPP}. 
\ei

Moreover, 
\bi
\item[$(iii)$] if $\lambda_1 < 0$ and $\delta \leq \delta_1$, there exists a positive stationary significant solution $p_{\delta}$ of \eqref{Eq_KPP} such that
 \be\label{p_delta}
 p_{\delta} \ge - \frac{\lambda_1 \phi_1}{\overline{\nu} (1+ \underline{\phi_1})}.
 \ee 
\item[$(iv)$] if $\lambda_1 < 0$ and $\delta > \delta_2$, the only possible positive bounded stationary solutions of \eqref{Eq_KPP} are remnant. 
\ei

{\bf Asymptotic behavior \cite{RC07}.}
Assume that the initial datum $y(0,x)$ for  Eq.~\eqref{Eq_KPP}  is taken to be  $p_0$, the unique nonnegative steady state of  Eq.~\eqref{Eq_KPP} when $\delta =0$ \cite{bhr1}. Then,  the solution $y(t,x)$ is non-increasing in $t$ and its asymptotic behavior is characterized as follows:
\bi
\item[$(v)$]  {\bf Population resilience:} If $\delta \le \delta^*$, the solution $y(t, \cdot)$ converges to $p_\delta$ uniformly in $\Omega$ as $t$ goes to infinite, where $p_\delta$ in \eqref{p_delta} is the maximal stationary significant solution of \eqref{Eq_KPP}.  

\item[$(vi)$] {\bf Population extinction:} If $\delta > \delta_2$, the solution $y(t, \cdot)$ converges uniformly in $\Omega$ to a remnant solution of \eqref{Eq_KPP}.  
\ei

\et
For a proof of the above results, see \cite[Theorem 2.6]{RC07} for $(i)$ and $(ii)$, \cite[Theorem 2.10]{RC07} for $(iii)$ and $(iv)$, and \cite[Theorem 2.11]{RC07} for $(v)$ and $(vi)$. In particular, $(iii)$ and $(iv)$ show that $\delta_1$ and $\delta_2$ defined in \eqref{Eq_deltas_def} provide easily computable bounds for the maximum allowable harvesting intensity $\delta^*$.  These bounds have been shown to be quite sharp for a broad range of habitat configurations \cite[Figure 4]{RC07}. 

Numerical simulations from \cite{RC07} based on the analytical estimates of $\delta^\ast$, strongly supported that environmental fragmentation of the habitat has a significant impact on the maximum sustainable yield associated with $\delta^\ast.$  Essentially, the more fragmented is the habitat, the more the population is susceptible to decline towards a remnant state under harvesting pressure.

\subsection{Optimal control of harvested population at the edge of extinction} \label{sect_KPP_control_formul}
We describe hereafter the optimal control problem we aim at solving for \eqref{Eq_KPP}.  The motivation is as follows. 
Assume that the population whose evolution is governed by \eqref{Eq_KPP}  is under threat of extinction due to e.g.~some illegal practices of harvesting. 
According to Theorem \ref{Thm_RC07_results} this situation is encountered when e.g.~$\delta > \delta_2$. Indeed, 
without additional external intervention, the population will eventually settle down to a remnant state and thus will be close to extinction.
Our goal is to prevent such a situation by exerting a certain control on the population by releasing new individuals from the same species coming e.g.~from a natural reserve, into the unprotected area $\Omega$, once  the population in that area drops below a certain warning threshold $P_c$.

The natural reserve is assumed to be limited in resources, and thus the rescue operation is itself under constraints. We assume that the reserve contains the same species that we are aiming at saving in the unprotected area and that the reserve's population is at an equilibrium, i.e.~that it has reached a steady state. Because no harvesting is exerted in the reserve, we assume that the population dynamics in the reserve is governed by  \eqref{Eq_KPP} with $\delta=0$,  namely the KPP equation in an heterogeneous environment \cite{skt}. In the general case, the domain $\widehat{\Omega}$ of the reserve  is different from that of the unprotected area, although it can share similar features such as a similar area and proportion of the population's habitat (in which $\mu >0$). More precise considerations about the reserve will be formulated in the numerical section. For the moment, we specify our goal in general terms useful to frame 
our optimal control problem. Thus, we are aiming at  transporting a fraction of the population from the  reserve  $\widehat{\Omega}$ to the unprotected area  $\Omega$. The question is to determine which fraction to use in order to save the population at extinction threat in $\Omega$ while not only avoiding to remove all the population in the reserve  but also not causing extinction within this same reserve.

{\bf Warning time $t=\tau$.} We assume also that the population in $\Omega$ before harvesting was exerted ($\delta=0$), had reached the stationary state $p_0$ and that the total population size can be monitored in $\Omega$.  Based on this monitoring, we assume  that the authorities declare the population at threat  when the following warning threshold $P_c$ is reached
\be\label{Def_Pc}
P_c= \beta \int_\Omega y(0,x) \d x,  \quad  0< \beta < 1, 
\ee 
where $y(0,x)=p_0$ is the initial density of the population at time $t=0$.

Recall that we assume $\delta > \delta_2$ in  \eqref{Eq_KPP} for the unprotected  area $\Omega$.  
As mentioned above, a situation for which $\delta$ would be greater than $\delta_2$ may arise due for instance to
illegal hunting practiced without respect of some quota and a lack of reliable monitoring by the authorities. 

In such a case, Theorem \ref{Thm_RC07_results}-($vi$)  guarantees that $y(t,x)$ converges uniformly towards a
remnant solution to \eqref{Eq_KPP} and because $y(t,x)$  is continuous in time and space,  we conclude that there exists a time instant, $t=\tau$, such that
 $\int_{\Omega} y(\tau,x) \d x = P_c$.  Note that we want to avoid to reach extinction so $P_c$ and thus $\beta$ in \eqref{Def_Pc} must be chosen such that
 \be\label{Pc_cond}
 P_c > \frac{\epsilon |\Omega|}{\underline{\phi_1}},
 \ee
 where $|\Omega|$ denotes the $d$-dimensional volume of $\Omega$. In other words one wishes to trigger an alert signal when the situation corresponds still to a significant population size according to \eqref{Def_significant_soln}.
The time instant $\tau$ serves us thus to trigger an alert signal from which one then starts to release new individuals from the natural reserve into the unprotected area. It will be called the warning time instant.

{\bf Extinction time $t=T$ and optimal control problem.} The goal  consists then of restoring, through an appropriate management plan exploiting the reserve,  the population to a safe dynamics away from extinction. For that purpose and based on the insights gained from Theorem \ref{Thm_RC07_results},  we aim at driving the population distribution to be as close as possible to a significant target population, chosen here to be the significant steady state $p_{\delta'}$ of \eqref{Eq_KPP} for some $\delta' \le \delta_1$ where $\delta_1$ is defined in \eqref{Eq_deltas_def}. Recall that  $\delta_1$ is a critical threshold that can be estimated in practice and that guarantees that the KPPH equation possesses a significant steady state (thus away from extinction) for any $\delta \leq \delta_1$.

 We are limited by time in our action. We want indeed to save the population from extinction that will occur at time $T$, corresponding to first time instant at which the maximum of the solution equals $\epsilon/\underline{\phi_1}$ if nothing is done\footnote{Because we assume $\delta >\delta_2$, we know that such a time instant exists as the population will become remnant; see Theorem \ref{Thm_RC07_results}-(vi).}.  Denoting  by $u(t,x)$ the number of new individuals brought from the reserve at time $t$ and position $x$, we consider then the cost functional
\be\label{Eq_J}
J(y,u) =\frac{1}{2} \int_\tau^T \big|y(t) - p_{\delta'}\big|^2_{L^2(\Omega)} \d t + \frac{\kappa}{2} \int_\tau^T \big|u(t)\big|^2_{L^2(\Omega)} \d t.
\ee
Here $y(t)$ denotes the solution in $L^2(\Omega)$  of the KPPH equation forced by $u(t,x)$ and such that $y(t)=y_0$ at $t=\tau$, and $\kappa>0$.
The first term in the right-hand side (RHS) of \eqref{Eq_J} is to enforce closeness to $p_{\delta'}$, while the second term is an energy-type term related to the effort of bringing the individuals from the reserve into the unprotected area. 

More precisely, we are aiming at addressing the following type of optimal control problems
\begin{subnumcases} {\label{Optimal_control_pbmain}}
\; \underset{u \in L^2(\tau,T; L^2(\Omega))}\min \; J(y,u) \label{min_J}\\
\; \textrm{where $y(t,x)$ is the solution of}\nonumber \\
\; \partial_t y = D \nabla^2 y + \mu(x) y - \nu(x) y^2 - \delta \rho_{\epsilon}(y) +  u(t,x), \; (t,x)\in [\tau, T] \times \Omega, \label{Eq_controlled_KPP1} \\
\; \frac{\partial y}{\partial {\boldsymbol{n}}}=0, \; (t,x)  \in [\tau, T]\times \partial \Omega, \label{Eq_controlled_KPP2}\\
\;\mbox{with $y(t,x)=y_0(x)$  at $t=\tau$.} \label{Eq_controlled_KPP3}
    \end{subnumcases}

As explained below, we restrict to controls lying within a subset of $L^2(\tau,T; L^2(\Omega))$, by introducing constraints on the control. To describe these constraints we enter into more specificities about the model setting. Throughout this article, we focus on a particular case of growth rate, $\mu$, defined to be
\be \label{Eq_mu}
\mu(x) =  m\chi^{}_\Lambda, \;\;  m>0,
\ee
where $\chi_\Lambda$ denotes the characteristic function of a (possibly disconnected) subdomain $\Lambda$  of  $\Omega$: 
\be \label{Eq_chi}
\chi^{}_\Lambda(x) = \begin{cases}
1, & \text{ if } x \in \Lambda,\\
0, & \text{otherwise}.
\end{cases}
\ee
Such a spatial dependence of the coefficient $\mu$ emphasizes that the population reproduces only when the individuals are within the subdomain $\Lambda$, while they may spread outside of $\Lambda$ due to the diffusion term $D\nabla^2 y$.

 The choice of the control $u$ in \eqref{Optimal_control_pbmain} is designed to act on this subdomain to enhance chances of natural growth.
For numerical applications we will consider piecewise constant controls in space while allowing for fluctuations in time. For that purpose, 
the domain $\Lambda$ is decomposed into mutually disjoints subdomains such that 
\be \label{Eq_Lambda1}
\Lambda = \Lambda_1 \cup \Lambda_2 \cup \cdots \cup \Lambda_K.
\ee
Then, by introducing  
\be \label{Eq_varphi}
\varphi_j = \frac{1}{\sqrt{|\Lambda_j|}}\chi^{}_{\Lambda_j}, \quad 1 \le j \le K,
\ee
where $|\Lambda_j|$ denotes the $d$-dimensional volume of $\Lambda_j$, we consider the following set of admissible controls: 
\bea \label{Eq_admissible}
\mathcal{U}_{\textrm{ad}}  \! = \! \bigg\{ \! (t,x)\mapsto \sum_{j=1}^K \Gamma_j (t) \varphi_j (x) \, :\, \Gamma_j \in L^2([\tau, T],\mathbb{R}),\; 0\leq \Gamma_j(t) \le C_j \text{ for a.e. $t\in [\tau, T]$} \! \bigg\},
\eea
where the $C_j$'s are positive constants which,  as explained below, are imposed by exploitation policy of the protected reserve. 
Denoting by $\mathcal{H} = L^2(\Omega)$, the optimal control problem \eqref{Optimal_control_pbmain} becomes, for this set of admissible controls,
\be \label{Eq_OPC}
\min J(y,u) \; \text{ subject to  $(y, u) \in L^2(\tau, T; \mathcal{H}) \times \mathcal{U}_{\textrm{ad}}$ solving \eqref{Eq_controlled_KPP1}-\eqref{Eq_controlled_KPP3}.}
\ee
Due to the definition of $\mathcal{U}_{\textrm{ad}}$, this optimal control problem is under constraints  on the control\footnote{Note that the set $\mathcal{U}_{\textrm{ad}}$ will be endowed with the induced topology from that of $L^2(\tau, T; L^2(\Omega))$.}. In the next section we will show that this optimal control problem can be addressed efficiently by means of finite-dimensional Galerkin approximations. 
For the moment we discuss how the constraints $C_j$  in \eqref{Eq_admissible} arise from  reserve management considerations.

\br\label{Rem_existenceOPC}
Note that for any given nonnegative continuous initial datum $y_0$ in \eqref{Eq_controlled_KPP3}, the existence of an optimal pair $(y(t;u^*),u^*)$ for the optimal control problem \eqref{Eq_OPC} can be ensured by using classical arguments such as found in \cite[Theorem 5.7]{Tro10} in combination with the maximum principle for \eqref{Eq_controlled_KPP1}-\eqref{Eq_controlled_KPP3}. 

Indeed, \eqref{Eq_controlled_KPP1}-\eqref{Eq_controlled_KPP3} can be put into the form \cite[Eq.~(5.8)]{Tro10} by taking $d(x,t,y) = - \mu(x) y + \nu(x) y^2 + \delta \rho_{\epsilon}(y)$ and $b \equiv 0$ therein. The assumptions collected into \cite[Assumption 5.6]{Tro10} are satisfied here except the requirement that $d_y(x,t, y) \ge 0$.  This latter condition is used in \cite{Tro10} to ensure the existence of solutions to \cite[Eq.~(5.8)]{Tro10} and also to guarantee the uniform boundedness of the solutions given by \cite[(5.10)]{Tro10}; see \cite[Assumption 5.4 and Theorem 5.5]{Tro10}\footnote{In \cite[Assumption 5.4]{Tro10}, the condition $d_y(x,t, y) \ge 0$ is stated as  $d$ being monotonically increasing with respect to $y$.}. However, for \eqref{Eq_controlled_KPP1}-\eqref{Eq_controlled_KPP3}, the existence of solutions and the desired uniform bound also hold as a consequence of the maximum principle here; see Appendix~\ref{Sect_cond_verification}. 
\er

\subsection{Control constraints from reserve management} \label{sect_KPP_control_formul2}
We assume that the population density in the reserve is at its steady state $\widehat{p}_0$, at the warning time $\tau$. The total population in the reserve is thus given by 
\be \label{Eq_Mtot}
M_{tot}=  \int_{\widehat{\Omega}} \widehat{p}_0(x) \d x.
\ee
A naive strategy, referred to as unsupervised below, consists of  assuming $\Gamma_j (t)$  to be set to a constant value $C_j$ to determine. Then the total (unsupervised) population, $M_T$, released  in the unprotected area $\Omega$ from the reserve, between times  $t=\tau$ and $t=T$, is 
\be
M_T= \Big(\sum_{j=1}^K C_j \sqrt{|\Lambda_j|} \, \Big) (T-\tau).
\ee
In practice, it is reasonable  to assume that an exploitation policy holds for the natural reserve independently of the rescue operation that is aimed to take place for $\Omega$.  As a result, only a fraction $\alpha$ ($0<\alpha<1$) of the total population  $M_{tot}$ in the reserve is allowed to be used for populating another area such as $\Omega$.

An optimal control planning $u^\ast(t,x)$ solving 
the optimal control problem \eqref{Eq_OPC} will be declared as efficient if it leads to an efficiency ratio 
\be \label{Eq_ratio_E}
E=\frac{\int_\tau^T \int_\Omega u^\ast  \d x \d t}{M_{tot}}< \frac{M_T}{M_{tot}}.
\ee
Assuming that the unsupervised strategy would exploit the allowable fraction $\alpha$,
we have
\be\label{Eq_C}
\Big(\sum_{j=1}^K C_j \sqrt{|\Lambda_j|} \, \Big) (T-\tau)=\alpha{M_{tot}}.
\ee
Finally, we will assume in what follows that the  $d$-dimensional volume $|\Lambda_j|$ are the same for different $j$'s, and that $C_j=C_{j'}=C$ for $j\neq j'$.
This way, Eq.~\eqref{Eq_C} gives the value of the constraints $C_j$ in $\mathcal{U}_{\textrm{ad}}$. By doing so we have thus an impartial 
way to compare the solution $u^\ast$ of the optimal control problem \eqref{Eq_OPC}, to that obtained from  a naive, unsupervised strategy consisting of taking $\Gamma_j(t)\equiv C$.
We have however to keep in mind that in this design {\it operatus} of the constraints,  the  set of admissible controls, $\mathcal{U}_{\textrm{ad}}$ depends on $\alpha$, since the $C_j$ do. As a result, any optimal solution $u^\ast$ to \eqref{Eq_OPC}  does also depend on the allowable fraction $\alpha$ from the reserve.   

From an exploitation viewpoint, we might ask what is the best fraction $\alpha$ to use in order to maximize the productivity of new individuals in the unprotected area, by the whole operation. An natural metric to answer this question consists of calculating the population ratio
\be \label{Eq_ratio_PR}
P_R=\frac{\int_{\Omega} y(T,u^\ast)  \d x}{\int_{\Omega} p_{\delta'} \d x}. 
\ee     
The ratio $P_R$ allows us to assess which fraction of the targeted total  population has been obtained (when driven by the optimal control planning $u^\ast$) when one reaches what would have been the extinction time $T$ if no action would have been taken. A natural question is then 
to analyze how the choice of $\alpha$ impacts $P_R$ and the efficiency ratio $E$, the goal being,  at this stage of the discussion, to minimize $E$ while maximizing $P_R$ as much as possible.  This later aspect will be discussed in Sec.~\ref{Sec_numerical_results} dealing with the numerical results.

 Also we want to make sure that the excision of population from the reserve does not lead to extinction there, neither. 
For that, we can rely here again on the theoretical understanding from the harvesting problem as recalled in Sec.~\ref{Sec_background2}. Let us assume that the removal of individuals from the reserve follows also a harvesting  law in $\widehat{\Omega}$ of the form  $\sum_{j=1}^K C_j \varphi_j(x) \rho_{\epsilon}(\hat{y})$, where  $\hat{y}$ denotes the population density in the reserve.
If one wants $u^\ast= \sum_{j=1}^K C_j \varphi_j$, one must have $\hat{y}(t,x)\geq  \epsilon/\underline{\phi_1}$ for $t$ in $[\tau,T]$ and $x$ in the reserve's domain\footnote{We assume that the reserve's habitat, $\widehat{\Lambda}$, is decomposed also into $K$ mutually disjoints subdomains each of same size than its corresponding counterpart $\Lambda_j$ in $\Omega$.} $\widehat{\Omega}$ to ensure $\hat{y}$ to be significant (see \eqref{Def_significant_soln}),  imposing thus a constraint on the population  on the reserve.

On the other hand, let $\hat{\delta}_1$ denotes the critical threshold ensuring resilience of the population  within the reserve this time; see Theorem \ref{Thm_RC07_results} again. A simple comparison argument shows that as soon as 
\be\label{Cond_constr_state}
\sum_{j=1}^K C_j \varphi_j(x)\leq \hat{\delta}_1 \bigg(\sum_{j=1}^K \varphi_j(x) \bigg),  \; x\in \widehat{\Omega},
\ee
then a significant steady state is ensured to exist, favoring resilience of the reserve's population. 

The inequality \eqref{Cond_constr_state} is equivalent to
\be\label{Eq_ctr_reserve}
\mbox{ For } 1\leq j \leq K, \;\; C_j \leq  \hat{\delta}_1, 
\ee
since the $\widehat{\Lambda}_j$ are mutually disjoint.  Thus the rescue operation for the unprotected area must comply with the constraint 
\eqref{Eq_ctr_reserve}, and also from what precedes, with  $\hat{y}(t,x)\geq \epsilon/\underline{\phi_1}$ for $t\in[\tau,T]$ and $x$ in the reserve's domain $\widehat{\Omega}$. In particular, \eqref{Eq_ctr_reserve} imposes a restriction on the population size $M_T$  excised from the reserve to satisfy $M_T\leq \hat{\delta}_1 (\sum_{j=1}^K \sqrt{|\widehat{\Lambda}_j|})(T-\tau).$ 
Such aspects regarding the reserve management will be also discussed in Sec.~\ref{Sec_numerical_results} below. For the moment we focus in the next section on the mathematical aspects of solving the optimal control problem \eqref{Eq_OPC} via Galerkin approximations, for a given set of admissible controls $\mathcal{U}_{\textrm{ad}}$ (and thus for a given allowable fraction  $\alpha$).

\section{Optimal control from Galerkin approximations} \label{Sect_Galerkin}

\subsection{Convergence results and error estimates about the optimal control}\label{Sec_31}
 In this section, we present convergence results and error estimates regarding  Galerkin approximations of the optimal control problem \eqref{Eq_OPC} constructed from eigenprojections. The synthesis of nearly optimal controls based on these Galerkin approximations is then provided in Section~\ref{Section_PMP}, following a standard Pontryagin Maximum Principle (PMP) approach.  

Denote the spectral elements of the eigenvalue problem \eqref{Eq_spectral_pb} by $\{(\lambda_j, e_j) \, : \, j \in \mathbb{N}\}$, where the eigenvalues $\lambda_j$'s are arranged in an increasing order, and the eigenfunctions $e_j$'s are normalized such that $\|e_j\|_{L^2(\Omega)} = 1$. Denote also the $N$-dimensional Galerkin approximation of the controlled state $y$ by
\be \label{Eq_uN}
y_N(t,x) = \sum_{i=1}^N \xi_i(t) e_i(x).
\ee
Recall that the control $u$ is written as 
\be \label{Eq_v}
u(t,x) = \sum_{i=1}^K \Gamma_i(t) \varphi_i(x).
\ee
Note that $y_N$ in \eqref{Eq_uN} depends on the initial datum and the control $u$ driving Eq.~\eqref{Eq_controlled_KPP1} (see \eqref{Eq_controled_Galerkin_v1} below). Throughout this article, the initial datum for the Galerkin approximation is taken to be $\Pi_N y_0$, where $\Pi_N$ denotes the projector onto the subspace spanned by the  first $N$ eigenmodes solving \eqref{Eq_spectral_pb}.  Dependence of $y_N$ on $\Pi_N y_0$ or $u$ will be made apparent depending on the context.

The Galerkin approximation of  \eqref{Eq_controlled_KPP1}-\eqref{Eq_controlled_KPP2} reads then:
\be \label{Eq_controled_Galerkin_v1}
\frac{\d \xi_i}{\d t} = -\lambda_i \xi_i  + \sum_{j,k=1}^N B_{jk}^i  \xi_j \xi_k - \delta  \Big \langle \rho_{\epsilon} \Big(\sum_{j=1}^N \xi_j e_j \Big), e_i \Big \rangle  +  \langle   u(t,\cdot), e_i \rangle, \quad t \in [\tau, T], 
\ee
where $i = 1, \ldots, N$, and
\be\label{Eq_Bijk}
B_{jk}^i =  - \langle \nu(\cdot) e_je_k, e_i \rangle = -  \int_{\Omega} \nu(x) e_i(x)e_j(x)e_k(x) \d x. 
\ee
Introducing a $K \times N$ matrix $\mathcal{M}$, whose elements are defined by
\be \label{Def_M}
\mathcal{M}_{ij} = \langle  \varphi_i, e_j \rangle = \int_{\Omega}  \varphi_i(x) e_j(x)  \d x,
\ee
we can rewrite \eqref{Eq_controled_Galerkin_v1} as 
\be \label{Eq_controled_Galerkin}
\frac{\d \xi_i}{\d t} = - \lambda_i \xi_i  + \sum_{j,k=1}^N B_{jk}^i  \xi_j \xi_k - \delta  \Big \langle \rho_{\epsilon} \Big(\sum_{j=1}^N \xi_j e_j \Big), e_i \Big \rangle  +  \sum_{j=1}^K \mathcal{M}_{ji}\, \Gamma_j(t), \quad t \in [\tau, T].
\ee

The cost functional $J_N$ associated with the $N$-dimensional Galerkin approximation \eqref{Eq_controled_Galerkin} is
\be \label{cost_JN_v1}
J_N(y_N,u) = \frac{1}{2} \int_\tau^T |y_N(t; \Pi_N y_0, u) - \Pi_N p_{\delta'}|_{L^2(\Omega)}^2 \d t   + \frac{\kappa}{2} \int_\tau^T |u(t)|_{L^2(\Omega)}^2 \d t.
\ee
By introducing
\be \label{Def_xi_Gamma}
\boldsymbol{\xi} = (\xi_1, \ldots, \xi_N)^{\textrm{T}}, \textrm{ and } \; \boldsymbol{\Gamma} = (\Gamma_1, \ldots, \Gamma_K)^{\textrm{T}}, 
\ee
we rewrite the cost functional $J_N$ given by \eqref{cost_JN_v1} as: 
\be \label{cost_JN}
J_N(\boldsymbol{\xi}, \boldsymbol{\Gamma}) = \frac{1}{2} \int_\tau^T |\boldsymbol{\xi}(t) - \boldsymbol{P}|^2 \d t   + \frac{\kappa}{2} \int_\tau^T |\boldsymbol{\Gamma}(t)|^2\d t, 
\ee
where $\boldsymbol{P} = (P_1, \ldots, P_N)^{\textrm{T}}$ with
\be
P_i = \langle p_{\delta'}, e_i \rangle, \quad i = 1, \ldots, N. 
\ee

In connection to the set of admissible controls $\mathcal{U}_{\textrm{ad}}$ defined by \eqref{Eq_admissible}, we introduce 
\bea \label{Eq_admissible_comp}
\mathcal{V}_{\mathrm{ad}}  \!=\! \bigg\{ \! t\mapsto (\Gamma_1 (t), \ldots, \Gamma_K (t))^{\textrm{T}}  \, :\, \Gamma_j \in L^2([\tau, T],\mathbb{R}),\; 0\le \Gamma_j(t) \le C_j \text{ for a.e. $t\in [\tau, T]$} \! \bigg\}.
\eea
The Galerkin approximation of the optimal control problem \eqref{Eq_OPC} is thus given by
\bea \label{Eq_OPC_Galerkin}
\min J_N(\boldsymbol{\xi}, \boldsymbol{\Gamma}) \; & \text{ defined in \eqref{cost_JN} subject to } (\boldsymbol{\xi}, \boldsymbol{\Gamma})  \in L^2(\tau,T; \mathbb{R}^N) \times \mathcal{V}_{\mathrm{ad}}  \text{ solving}   \\
&  \text{ the $N$-dimensional Galerkin system \eqref{Eq_controled_Galerkin} with $\boldsymbol{\xi}(\tau) = \Pi_N y_0$.}
\eea
We have then the following convergence results and error estimates, linking the optimal control problem \eqref{Eq_OPC} to its Galerkin approximation
\eqref{Eq_OPC_Galerkin}.

\bt \label{Thm_cve_KPP}

{\bf Error estimates about the optimal control.} Let us consider the optimal control problem \eqref{Eq_OPC} along with its Galerkin approximation \eqref{Eq_OPC_Galerkin}. Assume the initial datum $y_0$ in \eqref{Eq_controlled_KPP3} is strictly positive.  

Then, there exist $C >0$  such that the optimal control $u^\ast$ for \eqref{Eq_OPC} and the optimal control $u^\ast_{N}$ for the reduced problem \eqref{Eq_OPC_Galerkin} admit the following error estimate: 
\be\label{Est_contr_diff_KPP}
\int_\tau^T \big|u^\ast (t) - u^\ast_{N} (t)\big|_{L^2(\Omega)}^2 \d t \le C \left[\sqrt{T-\tau} + T-\tau  \right]
\bigg(\sum_{j=1}^{2}\bigg(\int_0^T \varepsilon_N^{j}(t,u^\ast,u_N^\ast)\mathrm{d} t\bigg)^{\frac{1}{2}}\bigg),
\ee
with 
\be\label{Eq_cv_epsilonN}
\underset{t\in[\tau,T]}\sup \varepsilon_N^{j}(t,u^\ast,u_N^\ast) \xrightarrow{N \rightarrow \infty} 0, \mbox{ for }\; j=1,2.
\ee

{\bf Convergence results.} Furthermore, the solution $y_N(t; \Pi_N y_{0}, u)$ of \eqref{Eq_controled_Galerkin} converges uniformly to the solution $y(t; y_0,u)$ of \eqref{Eq_controlled_KPP1}-\eqref{Eq_controlled_KPP3}: 
\be  \label{uniform_in_u_conv_KPP}
\lim_{N\rightarrow \infty}  \sup_{u\in \mathcal{U}_{\textrm{ad}}} \sup_{t \in [\tau, T]} \big|y_N(t; \Pi_N y_{0}, u) - y(t; y_0,u)\big|_{L^2(\Omega)} = 0.
\ee
\et

In Appendix \ref{Sect_cond_verification} we show that these convergence results and error estimates are the consequence of
general results of \cite{CKL17} dealing with the Galerkin approximations of nonlinear optimal control problems in Hilbert spaces and recalled in Appendix \ref{Sect_theory}, for the reader's convenience. 

More precisely, Theorem~\ref{Thm_cve_KPP} above is a consequence of  Theorem~\ref{Thm_controller_est} and Theorem~\ref{Thm:uniform_in_u_conv} presented in Appendix~\ref{Sect_theory}. Theorems~\ref{Thm_controller_est} and \ref{Thm:uniform_in_u_conv} deal with the optimal control of a broad class of nonlinear evolutionary equations containing the KPPH equation considered in this article. The required assumptions of Theorem~\ref{Thm_controller_est} and Theorem~\ref{Thm:uniform_in_u_conv} are verified in Appendix~\ref{Sect_cond_verification} for the optimal control  problems \eqref{Eq_OPC_Galerkin} and \eqref{Eq_OPC}.

In particular, the analysis shows that the constant $C$ in \eqref{Est_contr_diff_KPP} depends on the local Lipschitz constant (in $L^2(\Omega)$) of $\mathcal{G}(y)=\int_{\Omega} |y-p_{\delta'}|^2 \d x$  in a neighborhood of the origin, and also on the growth rate of $J$ in some appropriate norms of its arguments as given by \eqref{Eq_growth_onJ} in Appendix \ref{Sect_theory}.  The precise expression of $\varepsilon_N$ in \eqref{Eq_cv_epsilonN} is found in \eqref{Est_contr_diff}.

\subsection{PMP solution to the reduced optimal control problem \eqref{Eq_OPC_Galerkin}} \label{Section_PMP}

Thanks to Theorem \ref{Thm_cve_KPP}, it is sufficient to solve the reduced optimal control problem \eqref{Eq_OPC_Galerkin} to obtain nearly optimal solutions  as soon as $N$ is sufficiently large. We will see that for the KPPH model considered here, one can already obtain nearly optimal solutions with $N=1$ for certain specifications of the model such as the domain size, and choice of the model's parameters.

One can now use standard techniques from finite-dimensional optimal control theory to solve the above reduced optimal control problem \eqref{Eq_OPC_Galerkin}; see \cite{Bryson_al75,bonnard2003singular,Kirk12}.
Here, we follow an indirect approach by relying on the Pontryagin Maximum Principle (PMP); see e.g.~\cite{PBGM64,Kirk12}.

In that respect, we introduce the following Hamiltonian associated with the reduced optimal control problem \eqref{Eq_OPC_Galerkin}: 
\be \label{Def_HG}
H_G(\boldsymbol{\xi}, \boldsymbol{z}, \boldsymbol{\Gamma}) = \frac{1}{2} |\boldsymbol{\xi} - \boldsymbol{P}|^2 + \frac{\kappa}{2} |\boldsymbol{\Gamma}(t)|^2  + \boldsymbol{z}^{\textrm{T}} \boldsymbol{f}(\boldsymbol{\xi},\boldsymbol{\Gamma}),
\ee
where $\boldsymbol{z} = (z_1, \ldots, z_N)^{\textrm{T}}$ denotes the costate associated with the state $\boldsymbol{\xi}$, and $\boldsymbol{f}$ is the Galerkin vector field whose $i$-th component is given by the RHS of \eqref{Eq_controled_Galerkin}. That is, 
\be \label{Eq_f_comp}
f_i(\boldsymbol{\xi},\boldsymbol{\Gamma}) =  - \lambda_i \xi_i  + \sum_{j,k=1}^N B_{jk}^i  \xi_j \xi_k - \delta  \Big \langle \rho_{\epsilon} \Big(\sum_{j=1}^N \xi_j e_j \Big), e_i \Big \rangle  +  \sum_{j=1}^K \mathcal{M}_{ji}\, \Gamma_j(t), \quad i = 1, \ldots, N. 
\ee

Let
\bes
(\boldsymbol{\xi}^*, \boldsymbol{\Gamma}^*)  \in L^2(\tau,T; \mathbb{R}^N) \times \mathcal{V}_{\mathrm{ad}}  
\ees
be an optimal pair for the reduced optimal control problem \eqref{Eq_OPC_Galerkin}, and denote by $\boldsymbol{z}^*$ the costate associated with the state $\boldsymbol{\xi}^*$. It follows from the PMP that the triplet  $(\boldsymbol{\xi}^*, \boldsymbol{z}^*, \boldsymbol{\Gamma}^*)$  must satisfy the following constrained Hamiltonian system (see e.g.~\cite[Section 5.3]{Kirk12}):
 \begin{subequations}   \label{Pontryagin relation}
\begin{align}
& \begin{rcases}
\frac{\displaystyle \d \boldsymbol{\xi}^*}{\displaystyle \d t}  = \nabla_{\boldsymbol{z}} H_G(\boldsymbol{\xi}^*(t), \boldsymbol{z}^*(t), \boldsymbol{\Gamma}^*(t)) \\  
\frac{\displaystyle \d \boldsymbol{z}^*}{\displaystyle \d t} =  - \nabla_{\boldsymbol{\xi}} H_G(\boldsymbol{\xi}^*(t), \boldsymbol{z}^*(t), \boldsymbol{\Gamma}^*(t))  
\end{rcases}, \hspace{1.75em} (\text{Hamiltonian system for $(\boldsymbol{\xi}^*, \boldsymbol{z}^*)$}) \label{Pontryagin-a} \\
& H_G(\boldsymbol{\xi}^*(t), \boldsymbol{z}^*(t), \boldsymbol{\Gamma}^*(t))  \le H_G(\boldsymbol{\xi}^*(t), \boldsymbol{z}^*(t), \boldsymbol{\Gamma}(t)), \;\, \forall \; \boldsymbol{\Gamma} \in \mathcal{V}_{\mathrm{ad}}, \;\;  (\text{optimality condition})  \label{Pontryagin-b} \\
 & \boldsymbol{z}^*(T) =  0, \hspace{20em} (\text{terminal condition})    \label{Pontryagin-c}
\end{align}
\end{subequations}
where \eqref{Pontryagin-a} and \eqref{Pontryagin-b} hold for all $t$ in $(\tau, T)$.  Here $\nabla_x$ stands for the gradient operator in the $x$-direction.

\br
The optimality condition states that an optimal control must minimize 
the Hamiltonian $H_G$. In general it is a necessary condition and not a sufficient condition.  In the general case, there may be controls that satisfy  \eqref{Pontryagin-b} but that are not optimal controls. Yet the PMP may delineate a nonempty class of candidates. A triplet  $(\boldsymbol{\xi}^*, \boldsymbol{z}^*, \boldsymbol{\Gamma}^*)$ solution to  \eqref{Pontryagin relation} is called an extremal. Extremal solutions play an important role in optimal control theory; see \cite{bonnard2003singular}.  
Sufficient conditions for extremal to provide optimal controls can be found in \cite{hartl1995survey}.  For the problem at hand, since the cost functional
$J_N$ in \eqref{cost_JN} is quadratic in $\boldsymbol{\Gamma}$ and the dependence
on the control is linear for the control system \eqref{Eq_controled_Galerkin}, it is known that $\boldsymbol{\Gamma}^\ast$  obtained from such an extremal 
is actually the unique optimal control of the optimal control problem \eqref{Eq_OPC_Galerkin}; see e.g. \cite[Sec.~5.3]{Kirk12}.
\er

 From \eqref{Pontryagin relation}, we derive now, the explicit formula of the optimal control $\boldsymbol{\Gamma}^\ast$ based on the optimal costate $\boldsymbol{z}^\ast$.
To do so, we first remark that
\be
- \frac{\partial H_G(\boldsymbol{\xi}, \boldsymbol{z}, \boldsymbol{\Gamma})}{\partial \xi_i} = -(\xi_i - P_i) - \sum_{j=1}^N z_j \frac{\partial f_j(\boldsymbol{\xi},\boldsymbol{\Gamma})}{\partial \xi_i}.
\ee
Now since
\be
\frac{\partial f_j(\boldsymbol{\xi},\boldsymbol{\Gamma})}{\partial \xi_i} = - \lambda_i \delta_{ij} + \sum_{k=1}^N (B_{ik}^j + B_{ki}^j) \xi_k - \delta  \Big \langle \rho'_{\epsilon} \Big(\sum_{k=1}^N \xi_k e_k \Big) e_i, e_j \Big \rangle,
\ee
where $\delta_{ij}$ denotes the Kronecker delta, we get
\bea\label{Eq_g_comp}
  - \frac{\partial H_G(\boldsymbol{\xi}, \boldsymbol{z}, \boldsymbol{\Gamma})}{\partial \xi_i}  & =  -(\xi_i - P_i)  \\
  & - \sum_{j=1}^N z_j \Big(- \lambda_i \delta_{ij} + \sum_{k=1}^N (B_{ik}^j + B_{ki}^j) \xi_k - \delta  \Big \langle \rho'_{\epsilon} \Big(\sum_{k=1}^N \xi_k e_k \Big) e_i, e_j \Big \rangle\Big). 
\eea
Since \eqref{Eq_g_comp} does not depend on $\boldsymbol{\Gamma}$, we denote $- \partial H_G(\boldsymbol{\xi}, \boldsymbol{z}, \boldsymbol{\Gamma})/\partial \xi_i$ by a function $g_i(\boldsymbol{\xi},\boldsymbol{z})$. The costate equation 
becomes then $\d \boldsymbol{z}^*/ \d t=\boldsymbol{g}(\boldsymbol{\xi}^*,\boldsymbol{z}^*)$.
Thus the Hamiltonian system \eqref{Pontryagin-a} together with boundary condition \eqref{Pontryagin-c} becomes  

\bea \label{Eq_BVP_for_PMP}
& \frac{\displaystyle \d \boldsymbol{\xi}^*}{\displaystyle \d t} = \boldsymbol{f}(\boldsymbol{\xi}^*, \boldsymbol{\Gamma}^*), \quad t \in (\tau, T),\\
& \frac{\displaystyle \d \boldsymbol{z}^*}{\displaystyle \d t} = \boldsymbol{g}(\boldsymbol{\xi}^*,\boldsymbol{z}^*), \quad t \in (\tau, T),  \\
& \boldsymbol{\xi}^*(\tau)  = \Pi_N y_0,  \quad \boldsymbol{z}^*(T) = 0,
\eea
where the components of $\boldsymbol{f}$ and $\boldsymbol{g}$ are defined in the RHSs of \eqref{Eq_f_comp} and \eqref{Eq_g_comp}, respectively. The condition $\boldsymbol{\xi}^*(\tau)  = \Pi_N y_0$ is obtained by projecting the initial condition $y(\tau,x) = y_0$ in \eqref{Eq_controlled_KPP3}. 

We show next how $\boldsymbol{\Gamma}^*$ depends on the costate $\boldsymbol{z}^*$.
To do so, we first remark that by using $\mathcal{M}$ defined in \eqref{Def_M}, and the expressions of $H_G$ and the $f_i$ (see \eqref{Def_HG} and \eqref{Eq_f_comp}),  the optimality condition \eqref{Pontryagin-b} becomes 
\be \label{Eq_optimality2}
\frac{\kappa}{2} |\boldsymbol{\Gamma}^*(t)|^2  + (\boldsymbol{z}^*(t))^{\textrm{T}}\mathcal{M}^{\textrm{T}} \boldsymbol{\Gamma}^*(t)   \le \frac{\kappa}{2} |\boldsymbol{\Gamma}(t)|^2  + (\boldsymbol{z}^*(t))^{\textrm{T}} \mathcal{M}^{\textrm{T}} \boldsymbol{\Gamma}(t), \;\, \forall \; \boldsymbol{\Gamma} \in \mathcal{V}_{\mathrm{ad}}, \; t\in (\tau, T).
\ee
 Then by introducing
\be \label{Eq_def_w}
\boldsymbol{w}^* = \mathcal{M} \boldsymbol{z}^*,
\ee
 we get
\bea 
\frac{\kappa}{2} |\boldsymbol{\Gamma}(t)|^2  + (\boldsymbol{z}^*(t))^{\textrm{T}}  \mathcal{M}^{\textrm{T}} \boldsymbol{\Gamma}(t) & = \frac{\kappa}{2} |\boldsymbol{\Gamma}(t)|^2  + (\boldsymbol{w}^*(t))^{\textrm{T}} \boldsymbol{\Gamma}(t) \\
& = \frac{\kappa}{2} \sum_{j=1}^K \Big(\Gamma_j(t) + \frac{1}{\kappa}w^*_j(t)\Big)^2 -  \frac{1}{2\kappa} \sum_{j=1}^K w^*_j(t)^2.
\eea
As a consequence, the optimality condition \eqref{Eq_optimality2}  becomes
\be \label{Eq_optimality3}
\sum_{j=1}^K \Big(\Gamma^*_j(t) + \frac{1}{\kappa} w^*_j(t)\Big)^2 \le \sum_{j=1}^K \Big(\Gamma_j(t) + \frac{1}{\kappa} w^*_j(t)\Big)^2, \;\, \forall \; \boldsymbol{\Gamma} \in \mathcal{V}_{\mathrm{ad}}, \; t\in (\tau, T).
\ee
 Recalling the control constraints $0\le \Gamma_j \le C_j$  resulting from \eqref{Eq_admissible_comp},  one obtains that \eqref{Eq_optimality3} holds if and only if,  for all $j= 1, \ldots, K$,
\be  \label{Eq_optimality4} 
\Gamma^*_j(t)= h_j (\boldsymbol{w}^*(t)) \stackrel{\textrm{def}}{=} \begin{cases}
0, & \text{ if $-w^*_j(t) < 0$}, \\
- \frac{1}{\kappa} w^*_j(t), & \text{ if $0 \le - w^*_j(t) \le \kappa C_j$}, \\
C_j, & \text{ if $ -w^*_j(t) >  \kappa  C_j$},
\end{cases} \qquad  \; \forall \; t \in (\tau, T).
\ee

 Now by substituting the expression of $\boldsymbol{\Gamma}^*$  thus obtained and by using \eqref{Eq_def_w}, the problem \eqref{Eq_BVP_for_PMP} reduces to a boundary value problem (BVP) in the variables $\boldsymbol{\xi}^*$ and $\boldsymbol{z}^*$.
The synthesis of an optimal control $\boldsymbol{\Gamma}^*$ boils down thus to solving the following BVP:
\bea \label{Eq_BVP_for_PMP2}
& \frac{\displaystyle \d \boldsymbol{\xi}^*}{\displaystyle \d t} = \boldsymbol{f}(\boldsymbol{\xi}^*, \boldsymbol{h}(\mathcal{M \boldsymbol{z}^*})), \quad t \in (\tau, T),\\
& \frac{\displaystyle \d \boldsymbol{z}^*}{\displaystyle \d t} = \boldsymbol{g}(\boldsymbol{\xi}^*,\boldsymbol{z}^*), \quad t \in (\tau, T),  \\
& \boldsymbol{\xi}^*(\tau)  = \Pi_N y_0,  \quad \boldsymbol{z}^*(T) = 0,
\eea
where components of the function $\boldsymbol{h}$ are defined in \eqref{Eq_optimality4}.

 The optimal control $u^*_N$ to \eqref{Eq_OPC_Galerkin} in $\mathcal{U}_{\textrm{ad}}$ (see \eqref{Eq_admissible})
  is finally given
by
\be \label{Eq_uN_opt}
u^*_N(t,x) = \sum_{i=1}^K \Gamma^*_i(t) \varphi_i(x),
\ee
where $\boldsymbol{\Gamma}^*$ is given by  \eqref{Eq_optimality4} with  $\boldsymbol{w}^*$ therein obtained as $ \boldsymbol{w}^* = \mathcal{M} \boldsymbol{z}^*$  for $\boldsymbol{z}^*$ solving \eqref{Eq_BVP_for_PMP2}.

\br\label{Rem_unconstrained}
Note that, if the constraints $0\le \Gamma_j \le C_j$ on the admissible control are removed in $\mathcal{V}_{\mathrm{ad}}$ in \eqref{Eq_admissible_comp}, the optimality condition \eqref{Pontryagin-b} can be equivalently written as (see again \cite[Section 5.3]{Kirk12}):
\be \label{Pontryagin-b_unconstrained}
\nabla_{\boldsymbol{\Gamma}} H_G(\boldsymbol{\xi}^*, \boldsymbol{z}^*, \boldsymbol{\Gamma}^*)   = 0.
\ee 
The control law for the constrained case is given by \eqref{Eq_optimality4}. In contrast,  one obtains from \eqref{Pontryagin-b_unconstrained} the following control law for the unconstrained case: 
\be \label{Eq_optimality_unconstrained}
\boldsymbol{\Gamma}^*  = - \frac{1}{\kappa} \mathcal{M} \boldsymbol{z}^*,
\ee
with $\mathcal{M}$ given by \eqref{Def_M}.
Note that in practice \eqref{Eq_optimality4} reduces to \eqref{Eq_optimality_unconstrained} when the $C_j$ are sufficiently large. 
\er

\section{Nearly optimal controls from low-dimensional surrogates: Numerical results}\label{Sec_numerical_results}

\subsection{Numerical setup}\label{Sec_numsetup}
To illustrate our theoretical framework, and to simplify the reproducibility of the numerical results shown below, 
our numerical experiments take place in the case of one-dimensional environments. In that context, the unprotected area is a bounded, connected domain given by an interval, namely $\Omega=(0,\ell)$, with $\ell >0$. To simplify the analysis, we choose the domain $\widehat{\Omega}$ of the  reserve  to be also given by an interval of same length.  

We set $D = 1$,  $\ell = 1$, and $\nu = 0.2$  in the corresponding KPPH model \eqref{Eq_KPP}.
To account for the possible effects of heterogeneity of the habitat, we consider two cases of subdomain, $\Lambda$, appearing in the definition of $\mu$ given by \eqref{Eq_mu}: 
\bea \label{Eq_Lambda}
\text{Case I:} & \qquad \Lambda = [ 1/2, 1], \\
\text{Case II:} & \qquad \Lambda = [0, 1/4] \cup [3/4, 1].
\eea
Note that $|\Lambda| = 0.5$ for both cases, but Case I corresponds to an habitat more aggregated that Case II. 
In each case we set $m=2$ in \eqref{Eq_mu}.

Regarding the set $\mathcal{U}_{\textrm{ad}}$ of  admissible controls (see \eqref{Eq_admissible}), we divide $\Lambda$ into $K=8$ segments of equal length, such that each segment is of length $1/16$. 

Following Sec.~\ref{sect_KPP_control_formul}, we assume that for each case the population evolving in the unprotected area according to
\eqref{Eq_KPP} is under extinction threat because  $\delta = \delta_2(1+f)$, with $f$ chosen to be equal to $0.1$. 
One aims at restoring the population to a safe dynamics leading towards a significant steady state $p_{\delta'}$ appearing in the cost functional $J$ given by \eqref{Eq_J} for some $\delta'<\delta_1$.\footnote{We recall, that once $\delta' \leq  \delta_1$, there exists indeed  a significant steady state, $p_{\delta'}$, due to  Theorem~\ref{Thm_RC07_results}-($iii$).}
For that purpose we choose $\delta' = \delta_1(1-f)$, also with $f=0.1$.  Recall that $\delta_1$ and $\delta_2$ are critical harvesting thresholds defined in 
\eqref{Eq_deltas_def}; see Sec.~\ref{Sec_background2}.

The harvest function $\rho_{\epsilon}$ in \eqref{Eq_rho_conditions} is explicitly given here for $\epsilon = 0.05$ by 
\be \label{Eq_rho}
\rho_{\epsilon}(x) = \begin{cases}
1, & \text{if $x \ge \epsilon$}, \\
0.5 \sin(\pi (x-0.5 \epsilon)/\epsilon)+0.5, & \text{if $0 < x < \epsilon$}, \\
0, & \text{otherwise}. 
\end{cases}
\ee

Recall that the optimal control problem \eqref{Eq_OPC} is carried out over the time interval $[\tau, T]$. 
Following Sec.~\ref{sect_KPP_control_formul}, given a solution $y$ of the KPPH model \eqref{Eq_KPP} (emanating from $p_0$), the warning time $\tau$ is chosen according to $\int_0^{\ell} y(\tau,x) \d x =P_c$ with $P_c$ defined in \eqref{Def_Pc} for some $0<\beta <1$. 
We choose $\beta = 1/4$ in all the numerical experiments. In particular 
$P_c$ satisfies \eqref{Pc_cond} for Cases I and II.  In each case, the extinction time $T$  corresponds to the first time instant at which the maximum (over $(0,\ell)$) of the solution $y(t,x)$ equals $\epsilon/\underline{\phi_1}$. 
Note that $y_0$ used in the optimal control problem \eqref{Eq_OPC}   is $y(\tau,\cdot)$ (see \eqref{Eq_controlled_KPP3}).

\begin{figure}
\includegraphics[width=.95\textwidth,height=0.4\textwidth]{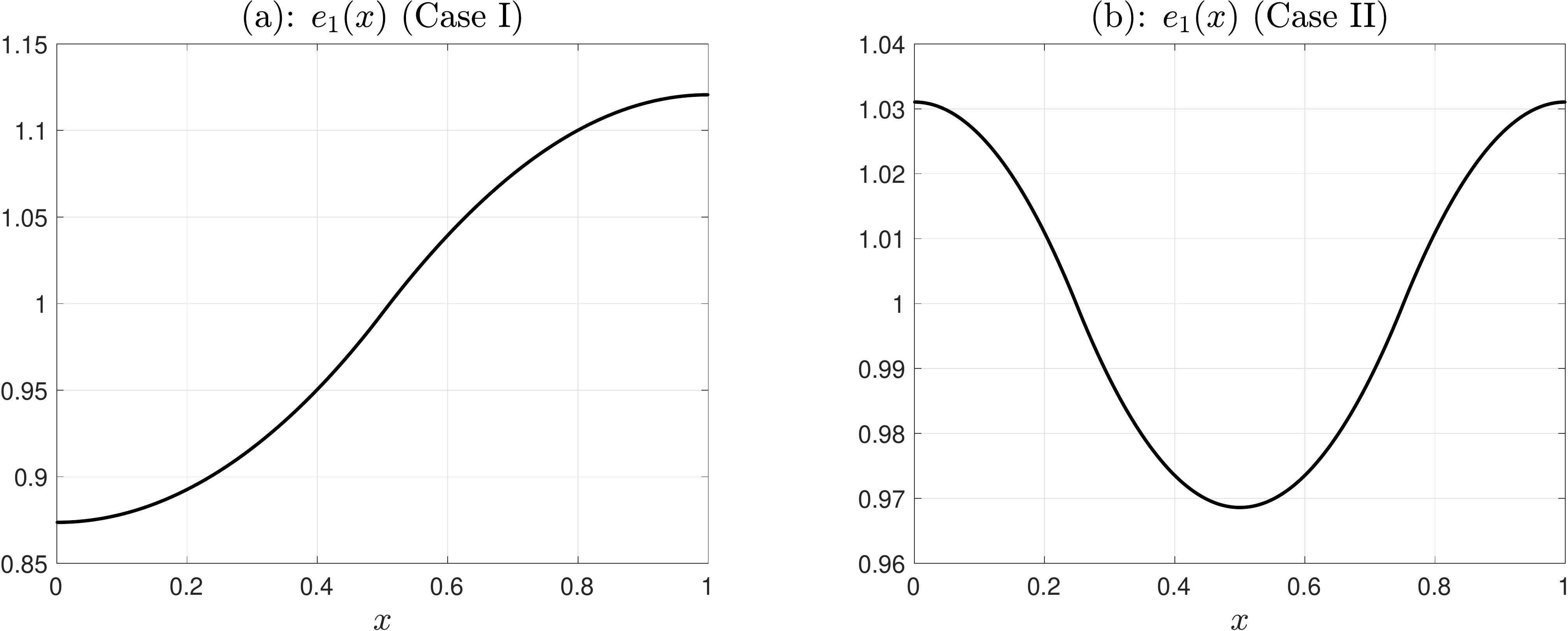}  
\caption{\footnotesize First eigenmode of the spectral problem \eqref{Eq_spectral_pb} for $\Omega=(0,1)$ shown here for $\Lambda=[ 1/2, 1]$ (Case I) and $ \Lambda = [0, 1/4] \cup [3/4, 1]$ (Case II).} 
\label{Fig_eigenmode}  
\end{figure}

The goal  is to drive, over the time window $[\tau, T]$,  the population distribution to be as close as possible to $p_{\delta'}$, while minimizing the cost functional given by \eqref{Eq_J} for controls $u$ lying in the admissible set $\mathcal{U}_{\textrm{ad}}$ given by \eqref{Eq_admissible}. The determination of the corresponding control constraints $C_j$ in \eqref{Eq_admissible} is discussed in Sec.~\ref{Choice_constraint} below. The rest of model's parameters for the two choices of subdomain $\Lambda$  given by \eqref{Eq_Lambda} (Cases I and II)  are listed in Table~\ref{Table_param_values},  rounded to the nearest thousandth. Note that $\delta_1$ and $\delta_2$ defined in \eqref{Eq_deltas_def} depend in each case on the first eigenmode  $e_1$ solving \eqref{Eq_spectral_pb} (associated with the first eigenvalue $\lambda_1$)\footnote{Note that $e_1$ is normalized in $L^2(\Omega)$ whereas $\phi_1$ appearing in \eqref{Eq_deltas_def}  is normalized in $L^\infty(\Omega)$. The two are related according to $\phi_1=e_1/\max(e_1)$.} shown here in Fig.~\ref{Fig_eigenmode}.

In what follows, the KPPH model is solved using the Matlab solver \texttt{pdepe} for $\Omega=(0,1)$ with $\delta x=10^{-2}$  and $ \delta t=10^{-2}$. 
The corresponding spectral problem \eqref{Eq_spectral_pb} is solved using finite difference with $\delta x =2^{-12}$. We use a small mesh grid size here to reach good accuracy given the discontinuity of the coefficient $\mu$. In both parameter cases, the dominant eigenvalue $\lambda_1$ is negative (see Table \ref{Table_param_values}), corresponding to an unstable eigenmode $\phi_1$. The rest of the spectrum corresponds here to stable eigenmodes.

\begin{table}[ht] 
\caption{\small Parameter values in \eqref{Optimal_control_pbmain}}
\label{Table_param_values}
\centering
\begin{tabular}{cccccccccc}
\toprule\noalign{\smallskip}
& $\lambda_1$ & $\delta_1$ & $\delta_2$ & $\delta'$ & $\delta$  & $\tau$  & $T$ \\  
\midrule\noalign{\smallskip}
Case I & -1.083 & 1.443 & 1.465 & 1.298 & 1.612 & 12.544 & 14.009 \\  
Case II & -1.021 & 1.301 & 1.302 & 1.171 & 1.432 & 14.023 & 15.582\\  
\noalign{\smallskip} \bottomrule 
\end{tabular}
\end{table}

\subsection{Choice of the control constraints}\label{Choice_constraint}
We assume that the reserve's domain size $|\widehat{\Omega}|$ is equal to $\ell$ and that the  population's habitat in the reserve $\widehat{\Lambda}$ satisfies $|\widehat{\Lambda}|=|\Lambda|.$  To simplify the analysis, we also assume that the  KPP model's parameters in the reserve $\widehat{\Omega}$ are the same as those of the KPPH equation in the unprotected area.
Under these working assumptions, the total population in the reserve is given by
\be \label{Eq_Mtot2}
M_{tot}=  \int_{\widehat{\Omega}} \widehat{p_0}(x) \d x,
\ee
where $\widehat{p_0}$ is the steady state of the KPP model (under Neumann condition) in $\widehat{\Omega}$.
In our numerical experiments, $M_{tot}= 5.34$ for Case I and $M_{tot}= 5.1$ for Case II.  

As pointed out in Sec.~\ref{sect_KPP_control_formul2}, the control planning $u(t,x)$ of transporting individuals from the natural reserve to the unprotected area is 
subordinated to the following reserve management factors:
\bi
\item[A)] The choice of the population fraction $\alpha$ allowed to be transported from the reserve.
\item[B)] The population displaced from the reserved must be such that the population distribution in the reserve satisfies $\hat{y}(t,x)\geq  \epsilon/\underline{\phi_1}$  for $t$ in $[\tau,T]$ and $x$ in $\widehat{\Omega}$.
\item[C)] For  $1\leq j \leq K,$ the constant $C_j$ must satisfy $C_j \leq  \hat{\delta}_1,$ where $\hat{\delta}_1$ is the critical harvesting threshold favoring population persistence for the reserve. 
\ei
 As explained in Sec.~\ref{sect_KPP_control_formul2}, the constraints $C_j$ in  \eqref{Eq_admissible} are determined according to \eqref{Eq_C}.
We choose $C_j=C$ for every $j$, i.e.~the constants are all the same across the subdomains $\Lambda_j$. Then  \eqref{Eq_C} allows us to find $C$ given by
\be\label{Eq_C2}
C =\frac{\alpha{M_{tot}}}{2(T-\tau)},
\ee
because $K=8$ and $|\Lambda_j|=1/16$. Note although the habitat dependence is not directly apparent in \eqref{Eq_C2}, the constant $C$ depends actually on $\Lambda$ through the coefficient $\mu$ impacting the dynamics and thus the extinction time $T$ and the warning time $\tau$.   
 Given our numerical setup, for the constant $C$ chosen according to \eqref{Eq_C2} then the requirement C) above is satisfied as long as $0<\alpha\leq \alpha_c $ with $\alpha_c=0.792$ for Case I  and  $\alpha_c=0.796$ for Case II. Requirement B) will be assessed a posteriori, after the optimal control problem is solved; see Sec.~\ref{Sec_numresults}.

\subsection{Control from effective reduced Galerkin systems} \label{Sec_4.3}
\begin{figure}
\includegraphics[width=.98\textwidth,height=0.4\textwidth]{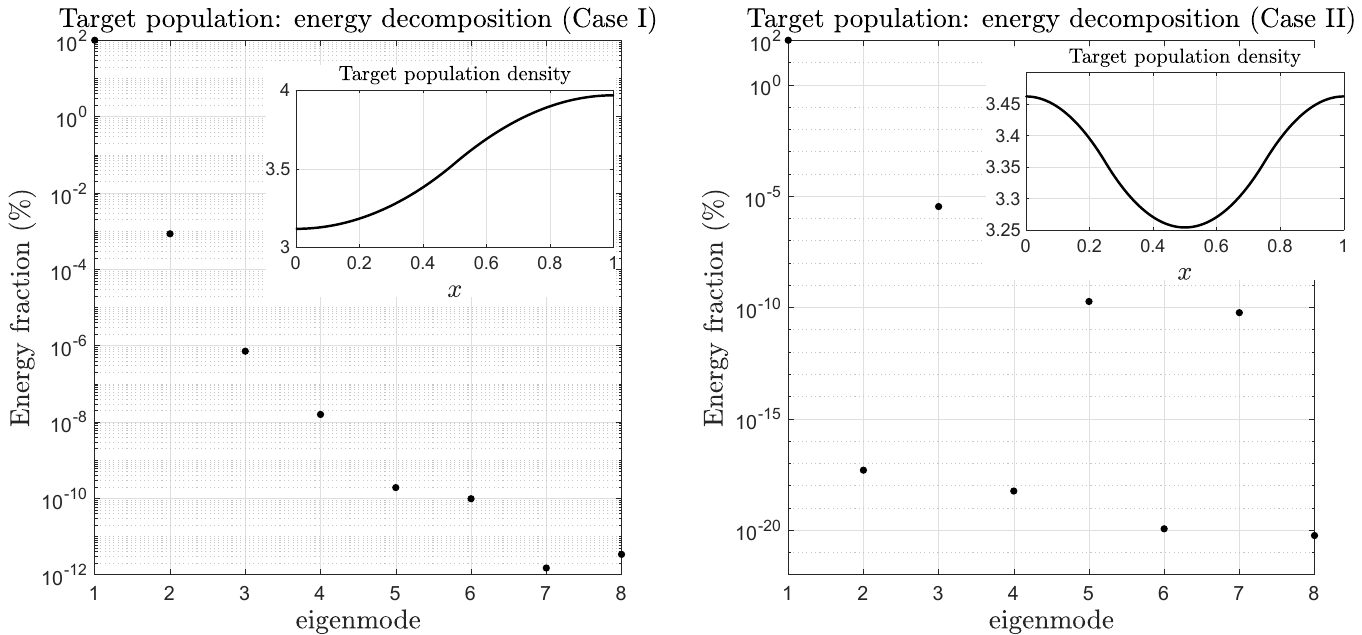}  
\caption{\footnotesize Energy decomposition of the target population density for $\Lambda=[ 1/2, 1]$ (Case I) and $ \Lambda = [0, 1/4] \cup [3/4, 1]$ (Case II) shown in semilog-scale. The target population density (shown in the upper insets) corresponds in each case to the significant steady states $p_{\delta'}$ of the KPPH equation \eqref{Eq_KPP} for $\delta=\delta'$, with $\delta'$ as indicated in Table \ref{Table_param_values}. These steady states are highly correlated with the corresponding first eigenmodes shown in Fig.~\ref{Fig_eigenmode}, explaining that these eigenmodes capture most of the energy contained in these steady states.} 
\label{Fig_energydecompostion}  
\end{figure}

To control the KPPH equation we approximate it by a Galerkin truncation of minimal dimension, namely we choose $N=1$. 
The reason is that for the model parameters considered here, the energy in the KPPH's solution is almost fully captured by the first eigenmode, the latter capturing more than 99.99\%  of the energy contained in the target population density $p_{\delta'}$ in each case; see Fig.~\ref{Fig_energydecompostion}. The insets of this figure show indeed that these target population densities are highly correlated with the corresponding first eigenmodes shown in Fig.~\ref{Fig_eigenmode}.

For $N=1$, the Galerkin approximation \eqref{Eq_controled_Galerkin} reduces to the following scalar logistic equation with harvesting and forcing terms:  
\be \label{Eq_controled_Galerkin_scalar}
\frac{\d X}{\d t} = X(a  -   bX)  - \delta g_\epsilon(X)+  \sum_{j=1}^K d_j \, \Gamma_j(t), \quad t \in [\tau, T], \; \; K=8,
\ee
where $a=-\lambda_1>0$, $b=-B_{11}^1$ (see \eqref{Eq_Bijk}) and the  $d_j=\mathcal{M}_{j1}$ ($N=1$) are constants given by \eqref{Def_M}, thus obtained here by the inner product between the eigenmode $e_1$ and the characteristic functions of the subdomains $\Lambda_j$. The harvesting term 
$g_\epsilon$ is given by
\be
g_\epsilon(X)= \begin{cases}
1& \text{ if }  X \geq \frac{\epsilon}{\min{e_1}},\\
\Big \langle \rho_{\epsilon} \Big(Xe_1 \Big), e_1 \Big \rangle & \text{if $X < \frac{\epsilon}{\min{e_1}}$}, 
\end{cases}
\ee
because the first eigenmode $e_1$ is positive. The harvest function $g_\epsilon$ satisfies the conditions in \eqref{Eq_rho_conditions}, and thus the Galerkin 
projection preserves the global structure of the original KPPH equation. 

Due to the high-energy content captured by the first-mode amplitude, $X(t)$, we expect that solving the optimal control \eqref{Eq_OPC_Galerkin}  associated with \eqref{Eq_controled_Galerkin_scalar} (i.e.~for $N=1$) should 
provide a nearly optimal solution for the original optimal control problem \eqref{Eq_OPC} with $\Omega=(0,\ell)$. The next section explores this intuition in more details.

Based on the analysis of Sec.~\ref{Section_PMP}, in what follows, the optimal control $u^*_N(t,x)$ (for $N=1$) is obtained according to \eqref{Eq_uN_opt}.  The optimal control $\Gamma_j^\ast$ therein is obtained  
from \eqref{Eq_optimality4}, namely 
$\Gamma^*_j(t)= h_j (\boldsymbol{w}^*(t))$, where  $\boldsymbol{w}^*=z^\ast \mathcal{M} $ with 
$z^\ast$ solving here the following BVP over $[\tau, T]$ (from \eqref{Eq_BVP_for_PMP2})
\bea \label{Eq_BVP_for_PMP3}
& \frac{\displaystyle \d X^*}{\displaystyle \d t} = X^*(a  -   bX^*)  - \delta g_\epsilon(X^*)+  \sum_{j=1}^K d_j h_j (z^\ast \mathcal{M} )\, ,\\
& \frac{\displaystyle \d z^*}{\displaystyle \d t} = -(X^\ast-\langle p_{\delta'}, e_1\rangle)-z^\ast(a-2b X^\ast)-  \delta k_\epsilon(X^\ast),  \\
& X^*(\tau)  = \langle y_0, e_1\rangle,  \quad z^*(T) = 0,
\eea
with
\be
 k_\epsilon(X^\ast)= \Big \langle \rho'_{\epsilon} \Big(X^\ast e_1 \Big) e_1, e_1 \Big \rangle.
\ee
Note that here $\mathcal{M}$ given by \eqref{Def_M} is a column vector and the costate $z^\ast$ is a scalar variable. 

In what follows, higher $N$-dimensional Galerkin truncation (with $N>1$) are also used  to control the KPPH equation. The use of such higher-dimensional  Galerkin  approximations is to benchmark the one-dimensional approximation.  The corresponding optimal control $u_N^\ast$
is, for each $N$, obtained  by following the PMP approach as described in Sec.~\ref{Section_PMP}.  
 
\begin{table}[htb] 
\caption{\small Coefficients of the 1D  Galerkin approximation \eqref{Eq_controled_Galerkin_scalar}}
\label{Table_coeffs_1D_Galerkin}
\centering
\begin{tabular}{cccccccccccc}
\toprule\noalign{\smallskip}
& $a$ & $b$ &  $\delta$  & $d_1$ & $d_2$ & $d_{3}$   & $d_4$ & $d_5$ & $d_6$ &  $d_7$ & $d_8$  \\ 
\midrule\noalign{\smallskip} 
Case I  & 1.083 &  0.202 & 1.612 & 0.252 & 0.259 & 0.265 & 0.270 & 0.274 & 0.277 & 0.279 & 0.280  \\
Case II  & 1.021 &  0.200   & 1.432 & 0.257 & 0.257 & 0.255 & 0.252 & 0.252 & 0.255 & 0.257 & 0.257 \\
\noalign{\smallskip} \bottomrule 
\end{tabular}
\end{table}

\subsection{First numerical results}\label{Sec_numresults}
We focus here on two particular choices of the allowable fraction $\alpha$ ($\alpha=0.1$ and $\alpha=0.35$) of individuals transported from the reserve's population to the unprotected area $\Omega,$ for each case of habitat $\Lambda$ in \eqref{Eq_Lambda}. Recall that the discrimination between the two cases of habitats considered in  Sec.~\ref{Sec_numsetup} is aimed at assessing the possible effects of fragmentation of the habitat on the rescue operation.
 
Below, when one writes  $J(y,u_N^*)$ (resp.~$J_N(\boldsymbol{\xi}^\ast,\boldsymbol{\Gamma}_N^*)$) it corresponds to the value of $J$ defined in \eqref{Eq_J} (resp.~$J_N$ defined in \eqref{cost_JN}) when the solution $y$ (resp.~$\boldsymbol{\xi}^\ast$) to the controlled KPPH model \eqref{Eq_controlled_KPP1}-\eqref{Eq_controlled_KPP3}  (to the $N$-dimensional Galerkin approximation \eqref{Eq_controled_Galerkin} of the controlled KPPH model) is driven by $u_N^\ast$ (resp.~$\boldsymbol{\Gamma}_N^\ast$).

For both cases of habitat (Cases I and II) and for both $\alpha$-values ($\alpha=0.1$ and $\alpha=0.35$), we have computed the relative error in the cost values $J$ between the optimal solutions obtained from the 1D and the 10D Galerkin approximations (see  Table \ref{Table_J}), namely
\be\label{rel_error} 
\texttt{error}=\frac{|J(y,u_1^*) - J(y,u_{10}^*)|}{J(y,u_{10}^*)}\times 100\%. 
 \ee
We found that in all cases, this relative error is almost negligible, bounded by $1.5 \times 10^{-4}\%$.

Furthermore, the convergence of the cost value $J_N(\boldsymbol{\xi}^\ast,\boldsymbol{\Gamma}_N^*)$ is achieved quickly as  the first few digits in the cost value has already converged when $N$ is increased; see Table \ref{Table_JN}. These numerical results confirm thus the intuition that due to the high-energy content carried out by the first eigenmode, a 1D Galerkin approximation is sufficient to obtain {\bf nearly optimal solutions}.   

\begin{table}[htb] 
\caption{\small Cost value $J(y,u_N^\ast)$}
\label{Table_J}
\centering
\begin{tabular}{cccccccccccc}
\toprule\noalign{\smallskip}
& {\footnotesize ($\alpha=0.1$, $N=1$)} & {\footnotesize ($\alpha=0.1$, $N=10$)} &  {\footnotesize ($\alpha=0.35$, $N=1$)}  & {\footnotesize ($\alpha=0.35$, $N=10$)}  \\
\midrule\noalign{\smallskip} 
Case I   & 4.168706 & 4.168700 & 3.203817 & 3.203754  \\  
Case II  & 3.952086 & 3.952086 & 2.996004 & 2.995997  \\  
\noalign{\smallskip} \bottomrule 
\end{tabular}
\end{table}

\begin{table}[htb] 
\caption{\small Cost value $J_N(\boldsymbol{\xi}^\ast,\boldsymbol{\Gamma}_N^\ast)$}
\label{Table_JN}
\centering
\begin{tabular}{cccccccccccc}
\toprule\noalign{\smallskip}
& $N=1$ & $N=2$ &  $N=3$  & $N=4$ & $N = 5$ & $N=10$  \\
\midrule\noalign{\smallskip} 
\, Case I, $\alpha=0.1$ & 4.159675 & 4.168644 & 4.168642 & 4.168651 & 4.168651 & 4.168651  \\ 
Case II, $\alpha=0.1$   & 3.953014 & 3.953014 & 3.953425 & 3.953425 & 3.953425 & 3.953426  \\
\midrule\noalign{\smallskip} 
\, Case I, $\alpha=0.35$ & 3.185567 & 3.203669 & 3.203668 & 3.203694 & 3.203693 & 3.203690  \\
Case II, $\alpha=0.35$ & 2.996189 & 2.996193 & 2.997253 & 2.997253 & 2.997253 & 2.997254  \\
\noalign{\smallskip} \bottomrule 
\end{tabular}
\end{table}

\begin{figure}
\includegraphics[width=\textwidth,height=0.55\textwidth]{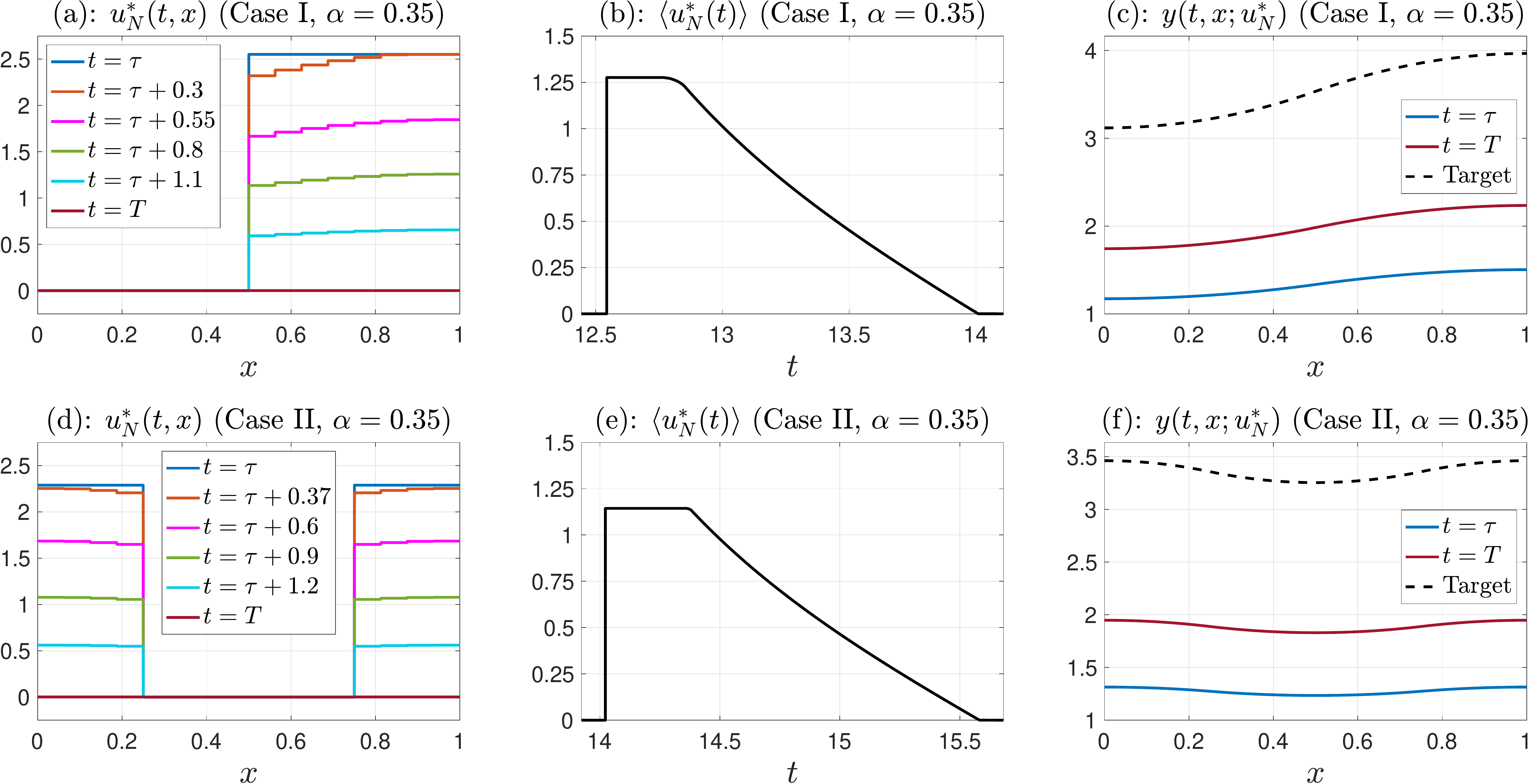}  
\caption{\footnotesize {\bf Panel (a)}: Optimal control $u_N^\ast(t,x)$ (Case I, $\alpha=0.35$) given  by \eqref{Eq_uN_opt} with $N=1$.  Here $\boldsymbol{\Gamma}_N^\ast$ is obtained  
from \eqref{Eq_optimality4}, where  $\boldsymbol{w}^*=z^\ast \mathcal{M} $ with 
$z^\ast$ solving here the BVP \eqref{Eq_BVP_for_PMP3}. {\bf Panel (b)}: Time-dependence of $\langle u_N^\ast(t)\rangle$ after space integration  over $(0,\ell)$.  The first time instant at which $\langle u_N^\ast(t)\rangle>0$ corresponds to $t=\tau$. The next time instant at which $\langle u_N^\ast(t)\rangle=0$ corresponds to $t=T$.  {\bf Panel (c)}: Controlled KPPH solution  $y(t,x;u_N^\ast)$ to  \eqref{Eq_controlled_KPP1}-\eqref{Eq_controlled_KPP3} when driven by $u_N^\ast$.
Panel (d) same as Panel (a), Panel (e) same as Panel (b), and Panel (f) same as Panel (c), for Case II.}
\label{Fig_alpha=dot35}  
\end{figure}

\begin{figure}
\includegraphics[width=\textwidth, height=0.55\textwidth]{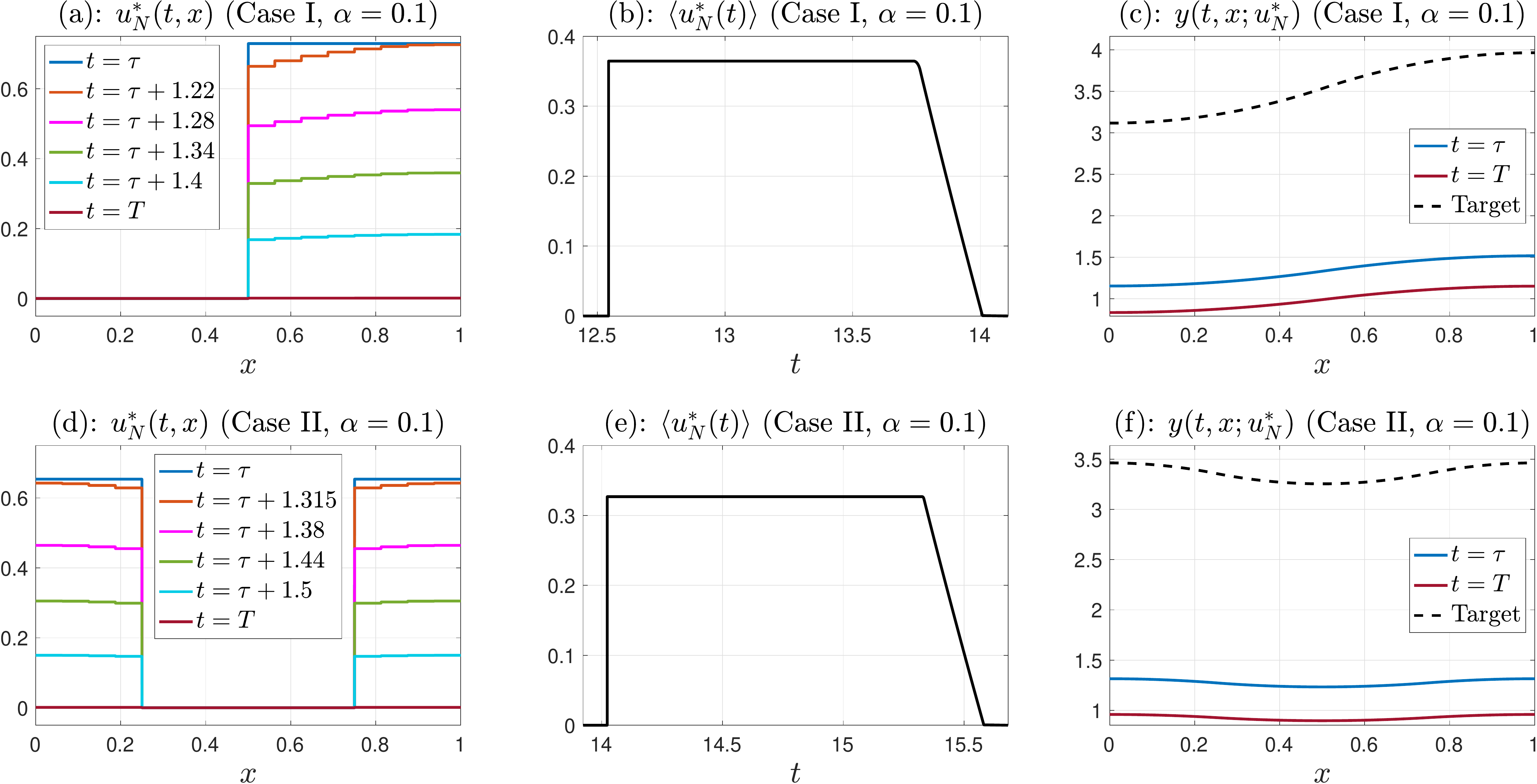} 
\caption{\footnotesize Same as Fig.~\ref{Fig_alpha=dot35} but for $\alpha=0.1$.} 
\label{Fig_alpha=dot1}  
\end{figure}

Based on these results,  we are in a comfortable position to discuss the numerical KPPH solution  $y(t,x;u_N^\ast)$ obtained when driven by the optimal control $u_N^\ast$ itself obtained according to \eqref{Eq_uN_opt} (for $N=1$) and the BVP \eqref{Eq_BVP_for_PMP3}. Figures \ref{Fig_alpha=dot35} and \ref{Fig_alpha=dot1} show these solutions in their respective panels (c) and (f) depending on the Case of habitat or value of $\alpha$ considered. The corresponding panels (a) (or (d)) show  the optimal control $u_N^\ast(t,x)$ (for $N=1$) for different time instants while panels (b) (or (e)) show its time evolution, after integration over $\Omega=(0,\ell)$, denoted by $\langle u_N^\ast(t)\rangle$.

Recall that given a value of $\alpha$, the constant $C$ is determined from \eqref{Eq_C2} which determines in turn the control constraints in the set $\mathcal{U}_{\textrm{ad}}$ of admissible controls, as the $C_j$ therein are chosen to be equal to $C$; see Sec.~\ref{Choice_constraint}. The impact of $\alpha$ is as follows. The smaller the allowable fraction $\alpha$ is, the smaller $C$ is and the more constrained is the optimal control problem \eqref{Eq_OPC}. This is visible on the temporal evolution of $\langle u_N^\ast(t)\rangle$ which is set to $C$ over a longer interval for $\alpha=0.1$ than for $\alpha=0.35$, before decreasing to zero; compare Panels (b) and (e) of Fig.~\ref{Fig_alpha=dot1} with those of Fig.~\ref{Fig_alpha=dot35}. We conduct a more detailed analysis on the dependence on $\alpha$ in Sec.~\ref{Sec_alpha_effects} below. For the moment, we discuss the effects of fragmentation of the habitat as the latter has been shown to play an important role in the population resilience to external perturbations \cite{RC07}.

\begin{table}[h] 
\centering
\caption{\small The efficiency ratio $E$ defined by \eqref{Eq_ratio_E}}  
\label{Table_E}
\begin{tabular}{ccc}
\toprule\noalign{\smallskip}
& $\alpha = 0.1$ &  $\alpha = 0.35$\\ 
\midrule\noalign{\smallskip}
Case I   & 0.0912 & 0.1920\\ 
Case II  & 0.0918 & 0.2002\\ 
\noalign{\smallskip} \bottomrule 
\end{tabular}
\end{table} 

The efficiency ratio $E$ defined in \eqref{Eq_ratio_E} serves us to compare the effects of different configurations of the habitat.
Recall that the smaller is this ratio the better it is in terms of exploitation of the reserve's population but not necessarily for saving the population from extinction in the unprotected area.  For that latter aspect one should have the population ratio $P_R$ (defined in \eqref{Eq_ratio_PR}) as large as possible while keeping in mind to respect the requirements A)-C) formulated in Sec.~\ref{Choice_constraint} in order to preserve the reserve's population from extinction.  

In terms of efficiency ratio, Table \ref{Table_E} indicates that, for a given allowable fraction $\alpha$ of the reserve, the 
habitat that is more homogeneous (Case I) is slightly more advantageous than the less homogeneous one (Case II). In terms of population ratio $P_R$ defined by \eqref{Eq_ratio_PR}, it is Case II that is slightly more advantageous, for instance $P_R=0.5633$ for Case II vs $P_R=0.5581$ for Case I, when $\alpha=0.35$. 
This means that for Case II, $56.33\%$ of the targeted population size is reached at time $t=T$, and $55.81\%$ for Case I.
Thus, the differences across the habitat fragmentation are somewhat minor.   To the opposite, the difference in terms of allowable fraction $\alpha$ of the reserve seem important as only of about $28\%$ (in both cases) of the targeted population size is reached at time $t=T$ when $\alpha=0.1$. The next section points out the main factors responsible of such a marked difference.

\subsection{Effects of the control constraints on the rescue operation}\label{Sec_alpha_effects}

When  $\alpha=0.1$, we observe  in Panels (c) and (f) of Fig.~\ref{Fig_alpha=dot1} that the final profile $y(T;u_N^\ast)$ of the controlled solution satisfies $y(T;u_N^\ast)<y_0$ ($y_0$ corresponds to the blue curve in Panels (c) and (f)), irrespectively of the degree of fragmentation of $\Lambda$. As time goes beyond $t=T$ and the harvesting intensity $\delta$ is set back to $\delta'$ in the KPPH equation while the control is abandoned\footnote{Assuming we get rid of  illegal harvesting by this time instant, while still allowing some harvesting respecting the quota $\delta' <\delta_1$.}, the population eventually settles down to a remnant state. The rescue operation has thus failed. 

On the contrary, we observe in Panels (c) and (f) of Fig.~\ref{Fig_alpha=dot35}  that for $\alpha=0.35$,  $y(T;u_N^\ast)>y_0$, and one gets in either Case I or II, closer than for $\alpha=0.1$ to the target population, i.e.~the significant steady state $p_{\delta'}$. As time goes beyond $t=T$ and $\delta$ is set to $\delta'$, the population converges towards this significant steady state (not shown, but see below). The rescue operation is a success.

What makes the difference so pronounced between $\alpha=0.1$ and $\alpha=0.35$? Remember that due to \eqref{Eq_C2} and our protocol for choosing our set $\mathcal{U}_{\textrm{ad}}$ of admissible controls, the allowable fraction $\alpha$ impacts the constraint $C$ arising in $\mathcal{U}_{\textrm{ad}}$. The smaller $\alpha$, the smaller $C$ and the more constrained becomes our control planning.   
At the same time it should be noted that when $\alpha$ is increased further towards $1$, the constraint $C$ becomes too large so that this constraint is never activated and the optimal control problem \eqref{Eq_OPC} behaves as  unconstrained.\footnote{It is easy to convince oneself of this fact by remarking that $\boldsymbol{\Gamma}^\ast$ given by \eqref{Eq_optimality4} reduces to \eqref{Eq_optimality_unconstrained} when the $C_j$ are sufficiently large; see Remark \ref{Rem_unconstrained}.} Thus there exists a critical  value of alpha, say $\overline{\alpha}$, above which the  control bounds are not activated in the set $\mathcal{U}_{\textrm{ad}}$ of admissible  controls. For Case I we found $\overline{\alpha}= 0.49$ and $\overline{\alpha}= 0.5$ for Case II. 
 
\begin{figure}
\centering
\includegraphics[width=.9\textwidth,height=0.35\textwidth]{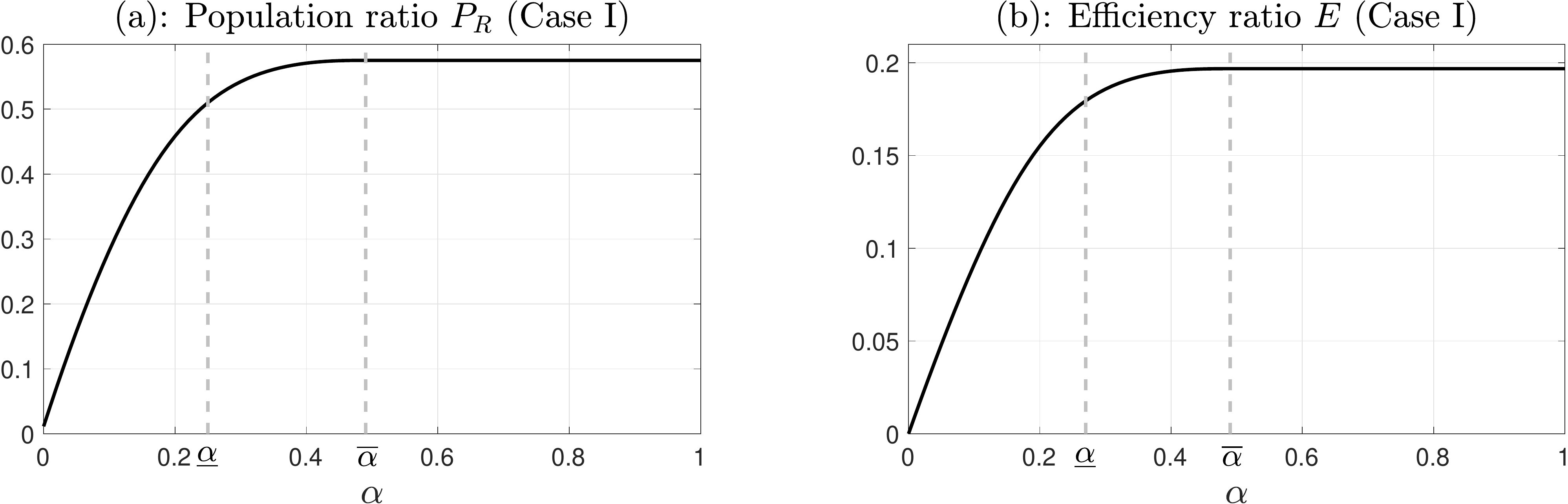}  
\caption{\footnotesize {\bf Panel (a)}: Population ratio $P_R$ (Case I). {\bf Panel (b):} Efficiency ratio $E$ (Case I). The value $\overline{\alpha}$ corresponds to the value of $\alpha$ above which the  control bounds are not activated.  The value $\underline{\alpha}$ corresponds to the value of $\alpha$ below which the population in $\Omega$ evolves towards a remnant state when the control is abandoned and $\delta$ is set to $\delta'$ after $t=T$; see \eqref{Cond_safe} below.}  
\label{Fig_PR}
\end{figure}

\begin{figure}
\centering
\includegraphics[width=.99\textwidth,height=0.45\textwidth]{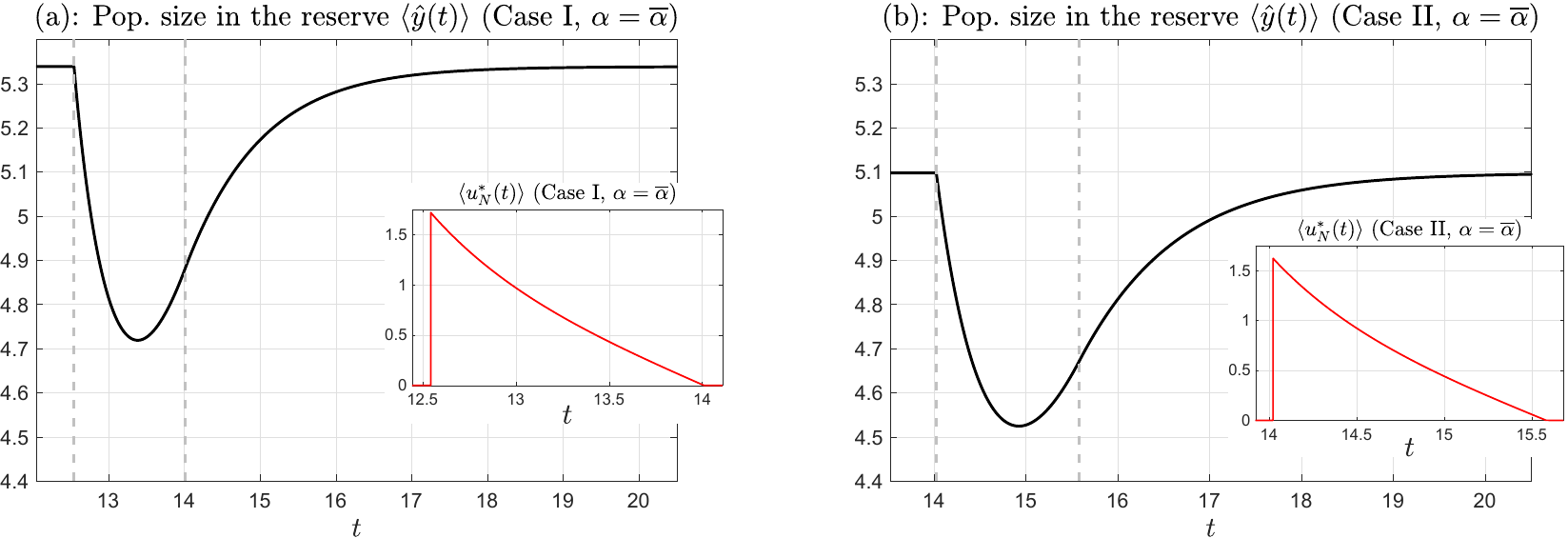} 
\caption{\footnotesize {\bf Panel (a)}: Population size $\langle \hat{y}(t)\rangle$ in the reserve  (Case I, $\alpha=\overline{\alpha}$). The bracket indicates integration over space.  The dynamics for the population in the natural reserve is governed by \eqref{Eq_KPP} forced by $-u^*_N$, corresponding to the amount of individual transported from the reserve to the unprotected area. The inset shows $\langle u_N^\ast (t) \rangle $ over a time window slightly larger than $[\tau,T]$. The first vertical dashed line indicates $t=\tau$, while the second one, indicates $t=T$. {\bf Panel (b)}: Same for Case II, $\alpha=\overline{\alpha}$.}  
\label{Fig_reserve}
\end{figure}

As a consequence, for  $\alpha>\overline{\alpha}$, only the same unconstrained optimal solution is obtained and the efficiency ratio $E$ saturates to a constant value.  
Panel (b) of Fig.~\ref{Fig_PR} shows the dependence of $E$ in terms of $\alpha$, and in particular its saturation for $\alpha>\overline{\alpha}$. Only Case I is shown here as $E$ for Case II behaves almost identically.  
The dependence on $\alpha$ of the population ratio $P_R$  is shown in  Panel (a) of Fig.~\ref{Fig_PR}. We observe that the $P_R$- and $E$-curves are highly correlated, both increasing and saturating for $\alpha>\overline{\alpha}$.  Thus, ``there is no free lunch,'' and one cannot minimize $E$ while maximizing 
$P_R$.  The choice of $\alpha$ is however a determining factor in the success of the rescue operation as pointed above. In all the cases, the excision of population to be transported from the reserve to the unprotected area can be absorbed by the reserve's population. Figure \ref{Fig_reserve} shows indeed that 
even for the more demanding ``excision pressure'' on the population's reserve (corresponding to the larger $E$  achieved for $\alpha=\overline{\alpha}$ and beyond), the population in the reserve recovers it original steady state, independently on the habitat's fragmentation.  

 The critical role of $\alpha$ on the survival of the population in the unprotected after application to optimal planning 
can gain great insights from the reduced system \eqref{Eq_controled_Galerkin_scalar} and the dynamical properties of the KPPH equation as recalled in Sec.~\ref{Sec_background2}. Let us assume that once the optimal control $u_N^\ast$ has been obtained from the BVP \eqref{Eq_BVP_for_PMP3} for a given $\alpha$, the harvesting intensity $\delta$ is set back to $\delta'$ in the KPPH equation, for $t>T$, and the control is abandoned. The idea is that by assuming that the illegal harvesting that was putting the population under extinction threat has been stopped starting from $t=T$, one aims at anticipating in terms of $\alpha$ (thus on the constraint put on the control) whether the rescue operation that would take place over the time window $[\tau,T]$, would be a success or not.

To do so, let us consider the following modification of \eqref{Eq_controled_Galerkin_scalar}
\be\label{Eq_proj_KPPH}
\frac{\displaystyle \d X}{\displaystyle \d t} = X(a  -   bX)  - \delta' g_\epsilon(X),
\ee
which is nothing else than the projection of the KPPH equation (with $\delta=\delta'$) onto the first eigenmode $e_1$. 
Let $\underline{X_\epsilon}$ denotes the smallest positive steady state of \eqref{Eq_proj_KPPH}. If $X(0)<\underline{X_\epsilon}$, then $X(t)$ converges 
to a remnant state whereas if $X(0)>\underline{X_\epsilon}$, it converges towards a significant steady state. 

\begin{figure}
\centering
\includegraphics[width=.99\textwidth,height=0.45\textwidth]{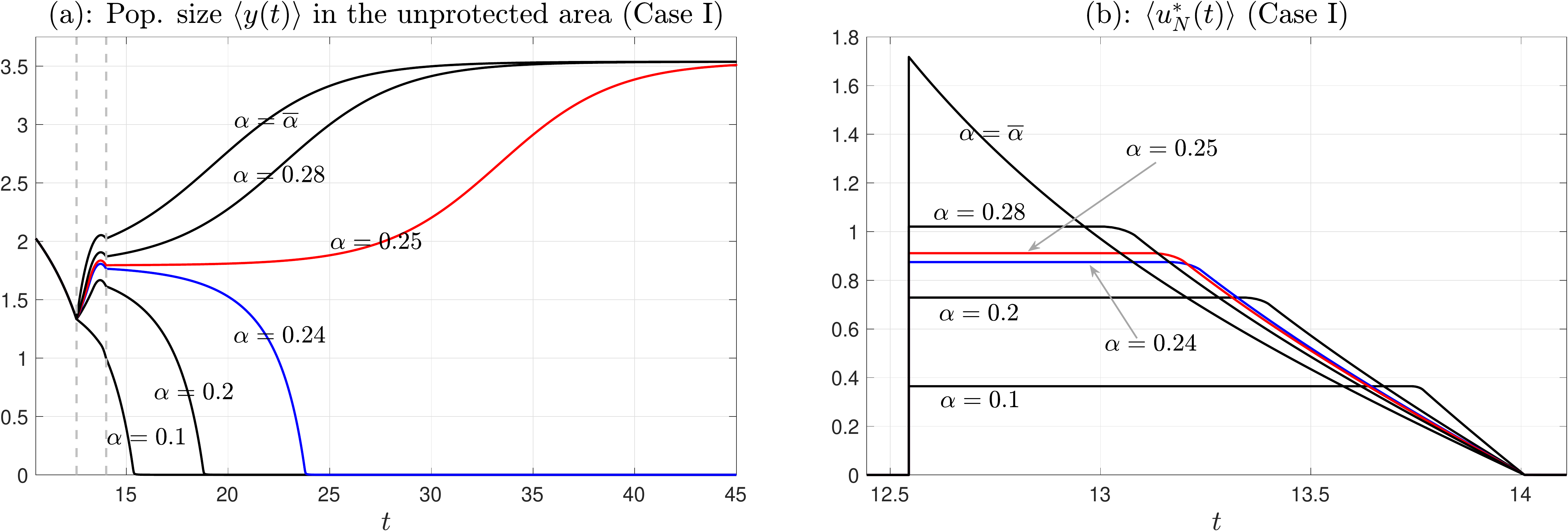}  
\caption{\footnotesize  {\bf Panel (a)}: Population size $\langle y(t)\rangle$ (after space integration) in the unprotected area $\Omega$ for Case I. 
The first vertical dashed line  indicates $t=\tau$, while the second one, indicates $t=T$.
For $t < \tau$, the dynamics is governed by \eqref{Eq_KPP} with $\delta$ set to its corresponding value given in Table~\ref{Table_param_values} (Case I). 
The initial datum at $t=0$ is taken to be $p_0$, the steady state of \eqref{Eq_KPP} for $\delta=0$. The dynamics over the time window $[\tau, T]$, is governed by \eqref{Eq_controlled_KPP1}--\eqref{Eq_controlled_KPP3} driven by $u_N^\ast (t)$ whose space integration is shown in Panel (b). For $t>T$, the dynamics is governed  by \eqref{Eq_KPP} but for $\delta = \delta'$ (Table~\ref{Table_param_values}) and with no control. {\bf Panel (b)}: Space integration, $\langle u_N^\ast (t) \rangle$,  of the optimal control $u_N^\ast$ obtained from  \eqref{Eq_uN_opt} (with $N=1$) and the BVP \eqref{Eq_BVP_for_PMP3}. Each curve is shown here as  $\alpha$ is varied according to the same values used for Panel (a) and for a time window slightly larger than $[\tau,T]$.}
\label{Fig_success}
\end{figure}

Marking the dependence of $u_N^\ast$ on $\alpha$, one defines $\underline{\alpha}$ to be the smallest value of $\alpha$ for which the following condition 
holds true
\be\label{Cond_safe}
\langle y(T;u_N^\ast(\alpha)),e_1\rangle> \underline{X_\epsilon}.
\ee
For Case I one finds $\underline{\alpha}=0.25$ using a mesh grid of size $0.01$ to discretize the range $[0,1]$ of $\alpha$-values. Because of the high-energy content carried by $e_1$, we expect that if \eqref{Cond_safe} is satisfied, then  the solution, $y(t;u_N^\ast(\alpha))$, to the KPPH equation for $\delta=\delta'$, converges to a significant state as $t$ tends to $\infty$.  
By showing the population size behavior (after space integration) $\langle y(t) \rangle$ as time evolves,  Panel (a) of Fig.~\ref{Fig_success} shows that this is exactly what happens for $\alpha\geq \underline{\alpha}$ while the population converges towards a remnant steady state for $\alpha<\underline{\alpha}$. One might thus wish to get $\alpha$ as close as possible to  $\underline{\alpha}$ from above in order to lower down $E$ while still ensuring success. However  one has to keep in mind that the same figure reveals that the boundary between 
success or failure of the rescue operation is very narrow, as $\alpha$ is getting too close to $\underline{\alpha}$. So in practice a hard constraint $C$ on the  individuals of the reserve should not correspond to an $\alpha$-value too close to  $\underline{\alpha}$, as $C$ may e.g.~not be exactly respected during the displacement operation. 

Panel (b) of Fig.~\ref{Fig_success} which shows the controls (applied over $[\tau,T]$) corresponding to the solutions shown in Panel (a), illustrates this statement.
As shown in this panel, the control leading to a significant survival (red curve) is indeed very close to that leading to extinction (blue curve). We are here in presence 
of an interesting phenomenon, namely, continuous dependence on the forcing may hold on finite-time intervals, but a high sensitivity in the system's response may take place in the asymptotic time.

We shall emphasize also that unlike what was reported for $\alpha=0.35$, one may have $y(T;u_N^\ast)>y_0$ while the population still evolves towards a remnant steady state, when $\delta$ is set to $\delta'$ after $t=T$ and the control is abandoned.    Figure \ref{Fig_success} shows for instance that the population size  in $\Omega$ on the rise from $y_0$ up to $t=T$ (due to control), drops eventually  to a remnant state in an asymptotic time for e.g.~$\alpha=0.2$ and $\alpha=0.24$. 

Finally, we stress that the numerical results reported are not limited to the particular numerical setup considered here. For instance if $D$ is further reduced, then more modes become unstable, say $p$, and the dimension of an efficient reduced system must be at least $p$.  Here $N=p=1$ for $D=1$. By setting
$D=0.1$ in Case II while keeping the other parameters as in Sec.~\ref{Sec_numsetup}, the two dominant eigenvalues of the spectral problem \eqref{Eq_spectral_pb} are negative corresponding to two unstable modes $e_1$ and $e_2$  capturing most of the energy of the target population density, unlike the energy decomposition shown in Fig.~\ref{Fig_energydecompostion}. By choosing $N=2$, the BVP \eqref{Eq_BVP_for_PMP3} becomes  \eqref{Eq_BVP_for_PMP2} for $N=2$. The corresponding efficient Galerkin approximation is then given by \eqref{Eq_controled_Galerkin_v1} for $N=2$ from which a critical $\alpha$ separating survival from extinction can also be determined from an analogue of  \eqref{Cond_safe} in which $\underline{X_\epsilon}$ (resp.~the projection onto $e_1$) is replaced by the steady state of smaller norm (resp.~the projection onto the space spanned by $e_1$ and $e_2$).  Such a remark about the inflation of the efficient reduced dimension holds as $D$ approaches zero. For instance for $D=0.01$ in Case II, there exists four unstable modes capturing most of the energy, constraining thus the efficient reduced dimension to be $N=4$. The approach presented here extends also to more realistic two-dimensional domains in which the fragmentation effects of the habitat may play a more important  role than for the one-dimensional domains; see \cite{RC07}.

\section*{Acknowledgments}
This work has been partially supported by the European Research Council (ERC) under the European Union's Horizon 2020 research and innovation program (grant agreement No.~810370 (MDC)), by the Office of Naval Research (ONR) Multidisciplinary University Research Initiative (MURI) grant N00014-20-1-2023 (MDC), and by the US National Science Foundation grant DMS-1616450 (HL).

\appendix 

\section{Galerkin approximations of nonlinear optimal control problems in Hilbert spaces: General convergence results} \label{Sect_theory}

In this appendix, we summarize from  \cite{CKL17} key convergence results and error estimates obtained 
for Galerkin approximations of nonlinear optimal control problems in Hilbert spaces. As we will see in Appendix \ref{Sect_cond_verification} below,
the material of this section allows us 
to link precisely the optimal control problem \eqref{Eq_OPC} to its Galerkin approximation
\eqref{Eq_OPC_Galerkin}, in terms of, both, error estimates about the optimal controls, and (strong) convergence about the controlled solutions as summarized in Theorem \ref{Thm_cve_KPP} in the Main Text.

In that respect, we tailor here the conditions of applications of the abstract results of \cite{CKL17} to the particular case of Galerkin systems from eigenprojections as considered in this article. The interested reader is referred to \cite{CKL17} for more general cases, and also for convergence results concerning the value functions. Applications to control in feedback form can be found in \cite{CKL18_DDE}.

In the following, we consider finite-dimensional approximations of the following initial-value problem (IVP):
\bea \label{Abstract_PDE}
\frac{\d y}{\d t} &= L y + F(y) + \mathfrak{C} (u(t)),  \quad t \in (0, T],\\
y(0) &= y_0,
\eea
where the unknown $y$ evolves in a separable Hilbert space $\cH$, $L\colon D(L) \subset \cH \rightarrow \cH$ is a linear operator with domain $D(L)$, $F\colon \cH \rightarrow \cH$ denotes the nonlinearity, and the initial datum $y_0$ belongs to $\cH$. The time-dependent forcing $u$ lives in a possibly different separable Hilbert space $V$; the (possibly nonlinear) mapping $\mathfrak{C}\colon V \rightarrow \cH$ is assumed to be such that $\mathfrak{C}(0)=0$.  Other assumptions regarding $\mathfrak{C}$ will be made below when needed. 

We assume the following conditions for the linear operator $L$. 
\bi
\item[{\bf (H0)}] {\it The linear operator $L: D(L) \subset \cH \rightarrow \cH$ is the  infinitesimal  generator of a $C_0$-semigroup of bounded linear operators $T(t)$ on $\cH$.} 

\item[{\bf (H1)}] {\it $L$ is self-adjoint with compact resolvent.} 
\ei
Recall that under assumption {\bf (H0)}, the domain $D(L)$ of $L$ is dense in $\cH$ and that $L$ is a closed operator;  see \cite[Cor.~2.5, p.~5]{Pazy83}.

For the admissible set of controls, we assume that 
\bi
\item[{\bf (H2)}] {\it Given some $q \ge 1$, the set of admissible controls $\mathcal{U}_{\textrm{ad}}$ is the subset of $L^q(0,T; V)$ constituted by measurable functions that  take values in  $U$, a bounded subset of the Hilbert space $V$.} 
\ei
In other words, 
\be\label{Eq_U_bounded}
\mathcal{U}_{\textrm{ad}}=\{f\in L^q(0,T; V)\;:  \; f(s) \in U \textrm{ for a.e. } s \in [0,T]\},\; \;q \ge 1.
\ee
The set $\mathcal{U}_{\textrm{ad}}$ will be endowed with the induced topology from that of $L^q(0,T; V)$.

Let $u$ be in $\mathcal{U}_{\textrm{ad}}$  given by \eqref{Eq_U_bounded}, a {\it mild solution} to \eqref{Abstract_PDE} over $[0,T]$ is a function $y$ in $C([0,T],\cH)$ such that
\be\label{Eq_mild}
y(t)=T(t) y_0 + \int_0^t T(t-s) F(y(s)) \d s + \int_0^t T(t-s) \mathfrak{C}(u(s)) \d s, \;\; t\in [0,T].
\ee
In what follows we will often denote by $t\mapsto y(t;y_0,u)$ a mild solution to \eqref{Abstract_PDE}.

Since $L$ is assumed to be self-adjoint with compact resolvent, it follows from spectral theory of self-adjoint compact operator  that the eigenfunctions of $L$ form an orthonormal basis of $\cH$; see e.g.~\cite{Brezis10}. We denote the eigenpairs of $L$ by $\{(\beta_k, e_k) \, : \, k \in \mathbb{N}\}$. For each $N \geq 1$, let $\cH_{N}$ be the $N$-dimensional subspace of $\cH$ spanned by the first $N$ eigenfunctions of $L$:
\be \label{Eq_HN}
\cH_{N}  = \mathrm{span}\{e_k \;:\; k = 1, \ldots, N\}. 
\ee
Denote also by $\Pi_N: \cH \rightarrow \cH_N$ the associated orthogonal  projector. Note that
\be\label{Eq_XN_in_domain}
\cH_N\subset D(L), \; \forall \, N\geq 1.
\ee

The corresponding Galerkin approximation of \eqref{Abstract_PDE} associated with $\cH_N$ is then given by:
\bea \label{ODE_Galerkin}
\frac{\d y_N}{\d t} &= L_N y_N + \Pi_N F(y_N) + \Pi_N \mathfrak{C} (u(t)), \; t\in [0,T],\\
y_N(0) &= \Pi_N y_0, \; \; y_0\in \cH,
\eea
where  
\be\label{Def_LN}
L_N = \Pi_N L \Pi_N : \cH \rightarrow \cH_N. 
\ee
In particular, the domain 
$D(L_N)$ of $L_N$ is $\cH$, because of \eqref{Eq_XN_in_domain}.

Throughout this  article, a mapping $f: \mathcal{W}_1 \rightarrow \mathcal{W}_2$ between two Banach spaces, $\mathcal{W}_1$ and $\mathcal{W}_2$, is said to be locally Lipschitz if for any ball $\mathfrak{B}_r \subset \mathcal{W}_1$ with radius $r>0$ centered at the origin,  there exists a constant $\Lip(f\vert{_{\mathfrak{B}_r}})>0$ such that 
\be  \label{Local_Lip_cond}
\|f(y_1) - f(y_2)\|_{\mathcal{W}_2} \le \Lip(f\vert{_{\mathfrak{B}_r}})\|y_1 - y_2\|_{\mathcal{W}_1}, \qquad \Forall y_1, y_2 \in \mathfrak{B}_r.
\ee

We will make also use of the following assumptions.
\bi
\item[{\bf (H3)}]  {\it The mapping $F: \cH \rightarrow \cH$ is locally Lipschitz in the sense given in \eqref{Local_Lip_cond}.}

\item[{\bf (H4)}] {\it Let $\mathcal{U}_{\textrm{ad}}$ be given by \eqref{Eq_U_bounded}.  For each $T>0$ and $(y_0,u)$ in $\mathcal{H}\times \mathcal{U}_{\textrm{ad}}$, the problem \eqref{Abstract_PDE} admits a unique mild solution $y(\cdot; y_0, u)$ in $C([0,T],\cH)$; and for each $N\geq 1$, its Galerkin approximation \eqref{ODE_Galerkin} admits a unique solution $y_N(\cdot; \Pi_N y_0, u)$ in $C([0,T],\cH)$. Moreover, there exists a constant $\mathcal{C}=\mathcal{C}(T,y_0)$ such that}
\begin{subequations} \label{Eq_y_uniform-in-u_bounds}
\begin{align}
& \|y(t; y_0, u)\|_{\cH} \le \mathcal{C}, && \Forall t\in[0,T], \; u \in  \mathcal{U}_{\textrm{ad}},  \label{Eq_y_uniform-in-u_bounds_a}\\
& \|y_N(t; \Pi_N y_0, u)\|_{\cH} \le \mathcal{C}, &&  \Forall t\in[0,T],\; N \in \mathbb{N}, \; u \in  \mathcal{U}_{\textrm{ad}}.  \label{Eq_y_uniform-in-u_bounds_b}
\end{align}
\end{subequations}

\ei

\br \label{Rmk:unif_bdd_on_soln}

Note that in applications, {\bf (H4)} is typically satisfied via {\it a priori} estimates. Indeed, let $u$ be in $\mathcal{U}_{\textrm{ad}}$  given by  \eqref{Eq_U_bounded}. Then uniform bounds such as in \eqref{Eq_y_uniform-in-u_bounds} are guaranteed if e.g.~an {\it a priori} estimate of the following type holds for the IVPs \eqref{Abstract_PDE} and \eqref{ODE_Galerkin}:
\be
\hspace{5ex}\underset{t \in [0,T]}\sup \| y(t;y_0,u)\|_{\cH}\leq \alpha ( \| y_0\|_{\cH}+ \|u\|_{L^{q}(0,T;V)}) +\beta,\;\; \;\; \alpha >0, \;\; \beta\geq 0.
\ee
See e.g.~\cite{Cazenave_al98,Tem97} for such a priori bounds for nonlinear partial differential equations.
Such bounds can also be derived for nonlinear systems of delay differential equations (DDEs); see \cite{CGLW15,CKL18_DDE}.
For the KPPH problem \eqref{Eq_controlled_KPP1}-\eqref{Eq_controlled_KPP2}, the desired estimate \eqref{Eq_y_uniform-in-u_bounds_a} is derived for nonnegative initial data by using the maximum principle together with energy estimates; see Appendix~\ref{Sect_cond_verification}. In particular as we will see, 
 Theorem \ref{Thm_controller_est} and Theorem \ref{Thm:uniform_in_u_conv} below,  both apply to this context when $y_0\geq 0$ in {\bf (H4)}.
\er

We introduce next the cost functional, $J\colon \cH \times \mathcal{U}_{\textrm{ad}}  \rightarrow \mathbb{R}^+$, associated with the IVP \eqref{Abstract_PDE}:
\be  \label{J_sec3}
J(y_0,u) = \int_0^T [\mathcal{G}(y(s; y_0, u)) + \mathcal{E}(u(s)) ] \, \d s,  \; y_0 \in \cH,
\ee 
where $\mathcal{G}: \cH \rightarrow \mathbb{R}^+$ and $\mathcal{E}: V \rightarrow \mathbb{R}^+$ are assumed to be {\it continuous}, and $\mathcal{G}$ is assumed to satisfy furthermore the condition:
\be\label{C1} 
\mathcal{G}  \mbox{ is  locally Lipschitz in the sense of } \eqref{Local_Lip_cond}.
\ee

The associated optimal control problem then writes 
\be  \label{P}  \tag {$\mathcal{P}$}
\begin{aligned}
\hspace{-.5ex}\min \, J(y_0,u)  \hspace{1ex}\text{ subject to } \hspace{1ex} (y, u) \in L^2(0,T; \cH) \times  \mathcal{U}_{\textrm{ad}} \text{ solves} ~\eqref{Abstract_PDE} \text{ with }  \; y(0)  = y_0 \in \cH.
\end{aligned}
\ee

The cost functional, $J_N\colon \cH_N \times \mathcal{U}_{\textrm{ad}}  \rightarrow \mathbb{R}^+$, associated with the Galerkin approximation  \eqref{ODE_Galerkin} is given by 
\be  \label{J_Galerkin}
J_N(\Pi_N y_0,u) = \int_0^T [\mathcal{G}(y_N(s; \Pi_N y_0, u)) + \mathcal{E}(u(s)) ] \, \d s,  \; y_0 \in \cH,
\ee 
and the corresponding optimal control problem reads: 
\be  \label{P_Galerkin}  \tag {$\mathcal{P}_N$}
\begin{aligned}
& \min \, J_N(\Pi_N y_0,u)  \quad \text{ subject to } \quad (y_N, u) \in L^2(0,T; \cH_N) \times  \mathcal{U}_{\textrm{ad}} \text{ solves} ~\eqref{ODE_Galerkin}\\
& \hspace{15em} \text{ with }  \; y_N(0)  = \Pi_N y_0 \in \cH_N.
\end{aligned}
\ee
Then, as a consequence of \cite[Corollary 2.14]{CKL17} we can deduce the following theorem about the error estimates between the full optimal control and the optimal control from  Galerkin approximations. 
\bt \label{Thm_controller_est}
Assume {\bf (H0)}--{\bf (H4)} and \eqref{C1} hold. Assume also that for each $y_0$ in $\cH$, both \eqref{P} and \eqref{P_Galerkin} admit an optimal control, denoted by  $u^*$  and $u^{*}_{N}$, respectively. Assume furthermore that there exists $\sigma >0$ such that the following local growth condition is satisfied  for the cost functional $J$ defined in \eqref{J_sec3}: 
\be\label{Eq_growth_onJ}
\sigma \|u^*  - v\|_{L^q(0,T; V)}^q \le J(y_0, v) - J(y_0, u^*), 
\ee
for all $v$ in some neighborhood $\mathcal{W} \subset \mathcal{U}_{\textrm{ad}}$ of $u^*$, with $\mathcal{U}_{\textrm{ad}}$ given by \eqref{Eq_U_bounded}. Assume finally that $u^*_N$  lies in  $\mathcal{W}$. 
Then there exists $\gamma >0$ such that
\bea\label{Est_contr_diff}
\|u^\ast - u^\ast_{N}\|_{L^q(0,T; V)}^q & \le \frac{1}{\sigma}\Lip(\mathcal{G}\vert_{\mathfrak{B}}) \left[\sqrt{T} + \gamma T  \right] \Bigl( \| \Pi_N^\perp y(\cdot; y_0, u^*)\|_{L^2(0,T; \cH)}\\ 
& \hspace{10em}+  2 \|  \Pi_N^\perp y (\cdot; y_0, u^{*}_{N})\|_{L^2(0,T; \cH)} \Bigr),
\eea
where $\mathfrak{B}$ denotes the ball in $\cH$ centered at the origin with radius $\mathcal{C}$, with $\mathcal{C}$ being the same as given in Assumption {\bf (H4)}, and $\Pi_N^\perp = \mathrm{Id}_{\cH} - \Pi_N$.

\et

As pointed out in \cite[Remark 2.13]{CKL17}, it is not clear a priori that 
\be \label{unif_vanishing}
\lim_{N\rightarrow \infty} \|\Pi_N^\perp y(\cdot; y_0, u^{\ast}_N)\|_{L^2(0,T; \cH)} = 0. 
\ee  
The reason relies on the dependence on $u^{\ast}_N$ of $\|\Pi_N^\perp y(\cdot; y_0, u^{\ast}_N)\|_{L^2(0,T; \cH)}$, where $u^{\ast}_N$ denotes the control synthesized from the $N$-dimensional Galerkin approximation. However, for the special case of Galerkin approximations constructed from eigenbasis, the convergence in \eqref{unif_vanishing} is guaranteed to hold under the assumptions {\bf (H0)}--{\bf (H4)} of Theorem \ref{Thm_controller_est} above; see \cite[Lemma 2.16]{CKL17}. As a result, when $N$ tends to infinity, 
$u^{\ast}_N$  converges to the optimal control $u^{\ast}$.

We have furthermore, the following uniform convergence result about the controlled solutions as a result of \cite[Theorem 2.6]{CKL17}.
\bt  \label{Thm:uniform_in_u_conv}
Assume that {\bf (H0)}--{\bf (H4)} hold and 
that the operator $\mathfrak{C}: V \rightarrow \cH$ is locally Lipschitz in the sense given in \eqref{Local_Lip_cond} (with $\mathfrak{C}(0) = 0$).
Assume furthermore that the set $U$ in Assumption {\bf (H2)} is a compact subset of $V$, with $q>1$ therein.  

Then, for any $(y_0,u)$ in $\cH\times \mathcal{U}_{\textrm{ad}}$, 
the solution $y_N$ of the Galerkin approximation \eqref{ODE_Galerkin} converges uniformly to the mild solution $y$ of \eqref{Abstract_PDE} in the sense that: 
\be  \label{uniform_in_u_conv_Goal}
\lim_{N\rightarrow \infty}  \sup_{u\in \mathcal{U}_{\textrm{ad}}} \sup_{t \in [0, T]} \|y_N(t; \Pi_N y_0, u) - y(t; y_0,u)\|_{\cH} = 0.
\ee
\et
Note that assumptions (A0), (A3) and (A6) in \cite[Theorem 2.6]{CKL17} correspond respectively to {\bf (H0)}, {\bf (H2)} and {\bf (H4)} assumed here, and assumption (A5) in \cite{CKL17} corresponds to {\bf (H2)} together with $U$ being compact.
Assumption (A7) required in \cite[Theorem 2.6]{CKL17} is actually here a consequence of {\bf (H0)}-{\bf (H4)} with  {\bf (H2)} assumed for $q>1$.

 Finally, we emphasize that assumptions (A1) and (A2) required by \cite[Theorem 2.6]{CKL17} follow from {\bf (H1)}. Indeed, as pointed out above, thanks to {\bf (H1)}, the eigenfunctions of $L$ form an orthonormal basis of $\cH$. Then, the operator $L_N$ defined by \eqref{Def_LN} as the eigen projection of $L$ onto $\cH_N$ given by \eqref{Eq_HN} clearly satisfies 
\be
\lim_{N\rightarrow \infty} \|L_N \phi - L\phi \|_{\cH} = 0, \quad \forall \, \phi \in D(L).
\ee 
Assumption (A2) in \cite{CKL17} is thus satisfied. Since $L$ is self-adjoint, it is also clear that the linear flow $e^{L_N t}: \cH_N \rightarrow \cH_N$ generated by $L_N$ satisfies 
\be
\|e^{L_N t}\| \le e^{\beta_1 t}, \quad N \in \mathbb{N},\, t \ge 0,
\ee
where $\beta_1$ is the largest eigenvalue of $L$. Defining the extension $T_N(t): \cH \rightarrow \cH$ of $e^{L_N t}$ to be   
\be
T_N(t) \phi = e^{L_N t} \Pi_N \phi + (\mathrm{Id}_{\cH} - \Pi_N) \phi, \quad \phi \in \cH,
\ee
Assumption (A1) in \cite{CKL17} follows then by choosing the parameters $M$ and $\omega$ therein to be $M=1$ and $\omega = \max\{\beta_1, 0\}$. 

\section{Convergence and error estimates results for the optimal control of the KPPH equation} \label{Sect_cond_verification}
We check here,  in the context of the optimal control of the KPPH equation, the assumptions of Theorem \ref{Thm:uniform_in_u_conv} and Theorem \ref{Thm_controller_est} about the convergence and 
error estimates results, respectively, allowing us in particular to deduce Theorem \ref{Thm_cve_KPP} of Sec.~\ref{Sec_31}. 

 Within this context, \eqref{Eq_OPC} and \eqref{Eq_OPC_Galerkin} play the role of \eqref{P} and \eqref{P_Galerkin}, respectively. The optimal problems \eqref{Eq_OPC} and \eqref{Eq_OPC_Galerkin}  are posed on the time interval $[\tau, T]$ but a simple change of variable 
$t' = t -\tau$ allows us to frame these as in Appendix~\ref{Sect_theory}, that is over the time interval $[0,T]$. We operate this shift below to ease 
the presentation. 

The verification of the required assumptions is organized in several steps.

\medskip
{\bf Step 1: Verifications of {\bf (H0)}, {\bf (H1)}, and {\bf (H3)}.} Let $\cH = L^2(\Omega)$. We first put the IVP \eqref{Eq_controlled_KPP1}-\eqref{Eq_controlled_KPP3} into the form \eqref{Abstract_PDE}. The corresponding operators $L \colon D(L) \rightarrow \cH$, $F\colon \cH \rightarrow \cH$ and $\mathfrak{C} \colon \cH \rightarrow \cH$ are then naturally defined as follows: 
\begin{equation}\label{Eq_def_operators}
\begin{aligned}
& L y = D \nabla^2 y + \mu(\cdot) y,  && y \in D(L) = H^2(\Omega) \cap \Big\{ y \in H^1(\Omega) \; \big\vert \; \frac{\partial y}{\partial \boldsymbol{n}}=  0 \Big \}, \\
& F(y) =  - \nu(\cdot) y^2 - \delta \rho_{\epsilon}(y), && y \in \cH, \\
& \mathfrak{C}(u) =   u,  && u \in \cH.
\end{aligned}
\end{equation}

It is standard that Assumptions {\bf (H0)} and {\bf (H1)} are satisfied for the elliptic operator $L$ defined in \eqref{Eq_def_operators}; see e.g.~\cite{Pazy83}. 
It can be checked that $F$ is locally Lipschitz as mapping from $\cH$ to $\cH$, in the sense of \eqref{Local_Lip_cond}. Thus Assumption {\bf (H3)}  is satisfied. It is also clear that $\mathfrak{C}$ defined above is Lipschitz on $\cH$ and satisfies furthermore that $\mathfrak{C}(0) = 0$.
 
\medskip
{\bf Step 2: Verification of {\bf (H2)}.}  Recall that the admissible set $\mathcal{U}_{\textrm{ad}}$ is defined by \eqref{Eq_admissible}.  Since $\varphi_j = \frac{1}{\sqrt{|\Lambda_j|}}\chi^{}_{\Lambda_j}$ (cf.~\eqref{Eq_varphi}), each admissible control $u$ in $\mathcal{U}_{\textrm{ad}}$ takes value in the following subset $U$ of $\cH$ at a given time instant: 
\be
U = \Big \{\psi \in \cH \; \big \vert \;  \psi(x) = \sum_{j = 1}^K \frac{b_j }{\sqrt{|\Lambda_j|}}\chi^{}_{\Lambda_j}(x), \; x\in \Omega, \; 0 \le b_j \le C_j \Big \}.
\ee
Thus $\mathcal{U}_{\textrm{ad}}$ in \eqref{Eq_admissible} can be rewritten as
\be\label{U_adB3}
\mathcal{U}_{\textrm{ad}}=L^2(0,T; U). 
\ee
 Taking the space $V$ in {\bf (H2)} to be $L^2(\Omega)$, it is clear that $U$ is a bounded set in $V$. Assumption {\bf (H2)} is thus verified with $q=2$ because of \eqref{U_adB3}. We show below that $U$ is furthermore a compact set in $V$, which is required in Theorem~\ref{Thm:uniform_in_u_conv}.
Every sequence $(\psi_{n})$ in $U$ has a convergent subsequence in $U$.  Indeed for each $1\leq j \leq K$, the corresponding $b_j^n$ has a convergent subsequence as taking value in the bounded and closed  subset $[0,C_j]$ of $\mathbb{R}$. Because we have a finite number of such convergent subsequences, one can always find a common extraction and thus $(\psi_{n})$ for which convergence holds towards an element of $U$.
Thus $U$ is a compact subset of $V = L^2(\Omega)$. 

\medskip
{\bf Step 3: Verification of {\bf (H4)}.} We first prove the uniform bound \eqref{Eq_y_uniform-in-u_bounds_a} for the solution $y$ to \eqref{Eq_controlled_KPP1}-\eqref{Eq_controlled_KPP3}. As it will appear below, this uniform bound can be derived as soon as 
one can ensure that $y\geq 0$ when $y_0\geq 0$. This property is a consequence of the maximum principle which holds for  
\eqref{Eq_controlled_KPP1}-\eqref{Eq_controlled_KPP3}, a known fact but of which we provide a proof of using the Stampacchia truncation method \cite{MS68}. The arguments are standard but are sketched below as they provide also useful insights to prove the required a priori bounds for the solution to the Galerkin approximation of \eqref{Eq_controlled_KPP1}-\eqref{Eq_controlled_KPP3}.
For this purpose, one defines the negative part $v^{-}$ of a measurable function $v: \Omega \rightarrow \mathbb{R}$ by
\be
v^{-}(x) = \min\{v(x), 0\}.
\ee
 A lemma due to Stampacchia ensures that if $v$ is in $H^1(\Omega)$, then $v^{-}$ lies also in $H^1(\Omega)$ and 
\be\label{Stamp_identity}
\nabla v^{-}=\chi_{\{v < 0\}} \nabla v. 
\ee
The so-called truncation method of Stampacchia allows us to obtain a powerful identity for a broad class of inhomgeneous heat problem 
\begin{subnumcases}{\label{Eqheat}}
\;\partial_t y = D \nabla^2 y +f, \quad  \textrm{ in } (0,T)\times \Omega, \label{Eqheat1}\\
\;\frac{\partial y}{\partial \boldsymbol{n}} =0, \quad  \textrm{ on }  (0,T)\times \partial \Omega, \\
\; y(0,x)=y_0 (x), \; x\in \Omega.
\end{subnumcases}
This identity is derived from the variational formulation of this problem, by using \eqref{Stamp_identity} with $v=y^{-}$. It gives
\be\label{Key_stamp}
\frac{1}{2} \big|y^{-}(t)\big|^2_{L^2(\Omega)} + D\int_0^t \big|\nabla y^{-}(s)\big|^2_{L^2(\Omega)} \d s=\int_0^t \int_{\Omega} f y^{-} \d x \d s+\frac{1}{2} \big|y^{-}(0)\big|^2_{L^2(\Omega)}.
\ee
From this identity the classical maximum principle can be deduce in the sense that  if $f\geq0$ and $y_0\geq 0$ a.e.~then $y\geq 0$ a.e.
Let us introduce
\be \label{Eq_N}
\mathcal{N}(x,y,u)  = \mu(x) y - \nu(x) y^2 - \delta \rho_{\epsilon}(y) +  u(t,x),
\ee
This result cannot be applied directly with $f=\mathcal{N}$ because $f$ does not have the good sign. However the identity \eqref{Key_stamp} is useful 
to conclude about the positivity of $y$ when $f=\widetilde{\mathcal{N}}$, with 
\be \label{Eq_N_modified}
\widetilde{\mathcal{N}}(x,y,u)  = \mu(x) y - \nu(x) y|y| - \delta \rho_{\epsilon}(y) +  u(t,x).
\ee
Indeed, first note that for $f=\widetilde{\mathcal{N}}$,
\be\label{Key_stamp2}
\int_0^t \int_{\Omega} f y^{-} \d x \d s \leq \|\mu\|_{L^\infty(\Omega)} \int_0^t \big| y^{-}(s)\big|^2_{L^2(\Omega)} \d s-\int_0^t \int_{\Omega} \nu(x) y|y| y^{-} \d x \d s
\ee
since $\rho_{\epsilon}(y)y^- = 0$, because $\rho_{\epsilon}(s) = 0$ for all $s \le 0$, and $u y^- \le 0$, because $u$ lies in $\mathcal{U}_{\textrm{ad}}$ and thus is nonnegative.
Then, since  
\be
\nu y |y| y^- = \nu |y| (y^-)^2 \ge 0,
\ee
we deduce from \eqref{Key_stamp2} that 
\be
\int_0^t \int_{\Omega} \widetilde{\mathcal{N}} y^{-} \d x \d s \leq  \|\mu\|_{L^\infty(\Omega)} \int_0^t \big| y^{-}(s)\big|^2_{L^2(\Omega)} \d s.
\ee
Assuming now that $y_0\geq 0$ a.e., we obtain from \eqref{Key_stamp} that 
\be
 \big|y^{-}(t)\big|^2_{L^2(\Omega)}\leq  2\|\mu\|_{L^\infty(\Omega)} \int_0^t \big| y^{-}(s)\big|^2_{L^2(\Omega)} \d s.
\ee 
Because $y^{-}(0)=(y_0)^-=0$, the Gronwall's lemma in its integral form allows us to conclude that $y^{-}=0$ a.e.~in $(0,T)\times \Omega$. 
Thus we have proved that $y\geq 0$ a.e.~if $y$ solves \eqref{Eqheat} with $y_0\geq 0$ and $f=\widetilde{\mathcal{N}}$ a.e.

From this we deduce that $|y|=y$ a.e.~$x$, $t$. Thus equation  \eqref{Eqheat1} with $f=\widetilde{\mathcal{N}}$  is the same as equation \eqref{Eq_controlled_KPP1} when 
$y_0\geq 0$ and in the end we have also found that $y$ is a positive solution to  \eqref{Eq_controlled_KPP1}-\eqref{Eq_controlled_KPP3} when $y_0\geq 0$. 

 On the other hand,  the inner product with $y$ on both sides of \eqref{Eq_controlled_KPP1} leads to
\bea 
\frac{1}{2} \frac{\d}{\d t} \big|y\big|^2_{L^2(\Omega)} & = \langle L y + F(y) + \mathfrak{C}(u(t)), y \rangle_{\cH} \\
& = - D  \big|\nabla y \big|^2_{L^2(\Omega)}  + \int_\Omega \big( \mu(x) y^2 - \nu(x) y^3 \big)  \d x - \delta \int_{\Omega}\rho_{\epsilon}(y) y \d x + \int_{\Omega} u y \d x
\eea
Thus because $y$ is a positive solution and $\rho_{\epsilon} \geq 0$, we obtain
\be\label{Eq_energy_est1}
\frac{1}{2} \frac{\d}{\d t} \big|y\big|^2_{L^2(\Omega)} \le \|\mu\|_{\infty} \big|y\big|^2_{L^2(\Omega)} +  C \int_{\Omega} y  \d x ,
\ee
with $C=\max \Big\{\frac{C_1}{\sqrt{|\Lambda_1|}}, \ldots, \frac{C_K}{\sqrt{|\Lambda_K|}} \Big\}$.

By remarking that by  H\"older's inequality, we have 
\bes 
\int_{\Omega} y  \d x\le \sqrt{|\Omega|} \big|y\big|_{L^2(\Omega)} \le \frac{1}{2}( |\Omega|+ \big|y\big|_{L^2(\Omega)}^2). 
\ees
Using this in \eqref{Eq_energy_est1}, we obtain then
\be
\frac{1}{2} \frac{\d}{\d t} \big|y\big|^2_{L^2(\Omega)}  \le c_1  \big|y\big|^2_{L^2(\Omega)}   + c_2,
\ee
where
\be
c_1 = \|\mu\|_{\infty} +  \frac{C}{2}, \quad c_2 = \frac{C }{2}|\Omega|. 
\ee
The uniform bound for $y$  in \eqref{Eq_y_uniform-in-u_bounds_a} follows then from Gronwall's inequality.

We turn now to the proof of the uniform bound \eqref{Eq_y_uniform-in-u_bounds_b}. First, we consider Galerkin approximation $\tilde{y}_N$ of 
a modified version of \eqref{Eq_controlled_KPP1}-\eqref{Eq_controlled_KPP2} in which the reaction term in  \eqref{Eq_controlled_KPP1} is replaced by $\widetilde{\mathcal{N}}$.  We assume the initial datum $\tilde{y}_0$ to be nonnegative.  
The corresponding $N$-dimensional Galerkin approximation $\tilde{y}_N$ satisfies thus   
\be \label{Eq_Galerkin_KPPH_modified}
\partial_t \tilde{y}_N = D \nabla^2 \tilde{y}_N + \Pi_N \widetilde{\mathcal{N}}(x,\tilde{y}_N,u).
\ee
From this equation one gets,
\bea 
\frac{1}{2} \frac{\d}{\d t} \big| \tilde{y}_N\big|^2_{L^2(\Omega)} & = - D  \big|\nabla  \tilde{y}_N \big|^2_{L^2(\Omega)}  +  \langle \Pi_N \widetilde{\mathcal{N}}(x, \tilde{y}_N,u),  \tilde{y}_N \rangle.
\eea
Since $\tilde{y}_N$ lies in $\cH_N$, we have 
\bea
\langle \Pi_N \widetilde{\mathcal{N}}(x,\tilde{y}_N,u), y_N \rangle & = \langle \widetilde{\mathcal{N}}(x,\tilde{y}_N,u), \tilde{y}_N \rangle \\
& = \int_\Omega \big( \mu(x) \tilde{y}_N^2 - \nu(x) |\tilde{y}_N|\tilde{y}_N^2 \big)  \d x - \delta \int_{\Omega} \rho_{\epsilon}(y_N) y_N \d x + \int_{\Omega} u \tilde{y}_N \d x \\
& \le  \|\mu\|_{\infty}  \int_\Omega \tilde{y}_N^2  \d x + (\delta + C) \int_{\Omega}  |\tilde{y}_N| \d x,
\eea
where $C$ is the same as in \eqref{Eq_energy_est1} and we have used the fact that $|\rho_{\epsilon}(y_N)|$ is bounded above by $1$; cf.~\eqref{Eq_rho_conditions}.

We infer thus that 
\bea 
\frac{1}{2} \frac{\d}{\d t} \big|\tilde{y}_N\big|^2_{L^2(\Omega)} \le  \|\mu\|_{\infty} \big|\tilde{y}_N\big|^2_{L^2(\Omega)} + (\delta + C) \int_{\Omega}  |\tilde{y}_N| \d x,
\eea
leading as for $y$ satisfying \eqref{Eq_energy_est1} to a uniform bound for $\tilde{y}_N$, and consequently \eqref{Eq_y_uniform-in-u_bounds_b} holds for $\tilde{y}_N$.
Due to Theorem \ref{Thm:uniform_in_u_conv}, $\tilde{y}_N$ satisfies the convergence property  \eqref{uniform_in_u_conv_Goal} towards the solution $\tilde{y}$ to
the modified IVP \eqref{Eq_controlled_KPP1}-\eqref{Eq_controlled_KPP3}  in which $\widetilde{\mathcal{N}}$ replaces the reaction term.  
Now we know that $\tilde{y}\geq 0$ since $\tilde{y}_0\geq 0$, and since $|\tilde{y}|=\tilde{y}$ in this case, we get that  $\tilde{y}$ solves actually
the original IVP \eqref{Eq_controlled_KPP1}-\eqref{Eq_controlled_KPP3} which we denote now by $y$ and $\tilde{y}_0$ by $y_0$.

Smoothing arguments ensure that the mild solution $y$ in $L^2(\Omega)$ to \eqref{Eq_controlled_KPP1}-\eqref{Eq_controlled_KPP3} is actually
continuous on $\Omega$. Thus because we have, for $\mathcal{U}_{\textrm{ad}}$ defined in \eqref{U_adB3} 
\be\label{Eq_totof}
\sup_{u\in \mathcal{U}_{\textrm{ad}}} \sup_{t \in [0, T]} |\tilde{y}_N(t; \Pi_N y_0, u) - y(t; y_0,u)|_{L^2(\Omega)} \underset{N\rightarrow \infty}\longrightarrow 0,
\ee
we have  $\tilde{y}_N > 0$ for $N$ sufficiently large (when $y_0>0$) and thus \eqref{Eq_Galerkin_KPPH_modified} reduces to the Galerkin approximation of the original IVP  \eqref{Eq_controlled_KPP1}-\eqref{Eq_controlled_KPP3}.

Finally, note that condition \eqref{C1} holds for the cost functional $J$ defined in \eqref{Eq_J}. The checking of \eqref{Eq_growth_onJ} is more involved 
but can be derived by adapting the proof of \cite[Theorem 5.3]{casas2017analysis} to our context.  The error estimates \eqref{Est_contr_diff} of Theorem \ref{Thm_controller_est} also hold with $q=2$ and $V=L^2(\Omega)$ and Theorem \ref{Thm_cve_KPP} of Sec.~\ref{Sec_31} is thus proved.

\bibliographystyle{amsalpha}
\providecommand{\bysame}{\leavevmode\hbox to3em{\hrulefill}\thinspace}
\providecommand{\MR}{\relax\ifhmode\unskip\space\fi MR }
\providecommand{\MRhref}[2]{%
  \href{http://www.ams.org/mathscinet-getitem?mr=#1}{#2}
}
\providecommand{\href}[2]{#2}

\end{document}